\documentclass[10pt]{scrartcl}
\usepackage{amsmath}
\usepackage{amssymb}
\usepackage{latexsym}
\usepackage{array}
\usepackage{yfonts}
\usepackage{fancyhdr}
\usepackage{epic, eepic}
\usepackage{graphicx,graphics}
\usepackage{color}
\usepackage{enumerate}
\DeclareMathAlphabet{\mathpzc}{OT1}{pzc}{m}{it}
\DeclareMathOperator{\vol}{vol}
\DeclareMathOperator{\eucl}{eucl}

\DeclareMathOperator{\spt}{spt}
\DeclareMathOperator{\grad}{grad}

\newtheorem{theorem}{Theorem}[section]
\newtheorem{defi}[theorem]{Definition}
\newtheorem{lemma}[theorem]{Lemma}

\newtheorem{cor}[theorem]{Corollary}
\newtheorem{rem}[theorem]{Remark}
\newtheorem{fig}{Figure}[section]

\newcommand{\R}{\mathbb{R}}

\newcommand{\N}{\mathbb{N}}

\newcommand{\beg}{\begin{equation}}
\newcommand{\en}{\end{equation}}
\newcommand{\beqn}{\begin{eqnarray}}
\newcommand{\eeqn}{\end{eqnarray}}
\newcommand{\beqnn}{\begin{eqnarray*}}
\newcommand{\eeqnn}{\end{eqnarray*}}

\newcommand{\K}{\mathrm{k}}
\numberwithin{equation}{section}

\author{}
\title{}
\date{}
\begin{document}
\maketitle
\vspace{-3.1cm}
\begin{center}{\LARGE Compactness of immersions with local Lipschitz representation}\end{center}
\vspace{7mm}
\begin{center} Patrick Breuning \footnote{P.\ Breuning was supported by the DFG-Forschergruppe
\emph{Nonlinear Partial Differential Equations: Theoretical and Numerical Analysis}.
The contents of this paper were part of the author's dissertation, which was written at
Universit\"{a}t Freiburg, Germany.}
\\ Institut f\"{u}r Mathematik der Goethe Universit\"{a}t Frankfurt am Main \\ Robert-Mayer-Stra{\ss}e 10,
D-60325 Frankfurt am Main, Germany \\email: breuning@math.uni-frankfurt.de \end{center} \vspace{1mm}
\begin{abstract}
\begin{center} {\textbf{Abstract}} \end{center}
\noindent We consider immersions admitting uniform representations
as a $\lambda$-Lipschitz graph. In codimension $1$, we show compactness for such
immersions for arbitrary fixed $\lambda<\infty$ and uniformly bounded volume. The same result is shown in arbitrary
codimension for $\lambda\leq \frac{1}{4}$.
\end{abstract}
%\vspace{1cm}
\begin{section}{Introduction}
In \cite{langer} J.\ Langer investigated compactness of immersed surfaces
in $\R^3$ admitting uniform bounds on the second fundamental form and the area of the surfaces.
For a given sequence $f^i:\Sigma^i\rightarrow \R^3$, there exist, after passing to a subsequence,
a limit surface
$f:\Sigma\rightarrow \R^3$ and diffeomorphisms $\phi^i:\Sigma\rightarrow \Sigma^i$, such that
$f^i\circ\phi^i$ converges in the $C^1$-topology to $f$.
In particular, up to diffeomorphism, there are only finitely many
manifolds admitting such an immersion.
The finiteness of topological types was generalized by K.\ Corlette in \cite{corlette} to immersions of arbitrary
dimension and codimension. Moreover, the compactness theorem was generalized by S.\ Delladio in
\cite{delladio} to hypersurfaces of arbitrary dimension. The general case, that
is compactness in arbitrary dimension and codimension, was proved by the author
in \cite{breuning1}. \\ \\
The proof strongly relies on a fundamental principle which we like to describe in the following. A
simple consequence of the implicit function theorem says that any immersion
can locally be written as the graph of a function $u:B_r\rightarrow\R^k$ over the affine tangent space. Moreover,
for a given $\lambda>0$ we can choose $r>0$ small enough such that $\|Du\|_{C^0(B_r)}\leq \lambda$. If
this is possible at any point of the immersion with the same radius $r$, we call $f$ an $(r,\lambda)$-immersion. \\ \\
Using the Sobolev embedding it can be shown that a uniform $L^p$-bound for the second fundamental form with $p$ greater
than the dimension implies that for any $\lambda>0$ there is an $r>0$ such that every immersion is an
$(r,\lambda)$-immersion. \\ \\
Inspired by this result, it is a natural generalization to investigate compactness properties also
for $(r,\lambda)$-immersions with fixed $r$ and $\lambda$; this is the topic of the present paper.
In the proof of the theorem of Langer it is essential that
$\lambda$ can be chosen very small. Then, using the local graph representation
over $B_r$, all immersions are close to each other and nearly flat.
These properties are used repeatedly, for example for the construction of the
diffeomorphism $\phi^i$. \\ \\
Here, we would first like to show compactness of $(r,\lambda)$-immersions in codimension $1$ for any fixed $\lambda$.
We do not require any smallness assumption for $\lambda$.
Moreover, we do not only consider immersions with graph representations over the affine tangent space, but also
over other appropriately chosen $m$-spaces.
Let $\mathfrak{F}^1(r,\lambda)$ be the set of $C^1$-immersions
$f:M^m\rightarrow \R^{m+1}$ with $0\in f(M)$, which may locally be written
over an $m$-space as the graph of a $\lambda$-Lipschitz function $u:B_r\rightarrow \R$
(the precise definitions of all notations used in this paper are given
in Section 2).
Here all manifolds are assumed to be compact. Moreover,
let $\mathfrak{F}^1_{^\mathcal{V}}(r,\lambda)$ be the set of
immersions in $\mathfrak{F}^1(r,\lambda)$ with $\vol(M)\leq \mathcal{V}$.
Similarly, we define the set $\mathfrak{F}^0(r,\lambda)$ by replacing
$C^1$-immersions in $\mathfrak{F}^1(r,\lambda)$ by Lipschitz functions.
We obtain the following compactness result: \\
\begin{theorem} \label{compactness2} (Compactness of
$(r,\lambda)$-immersions in codimension one) \\ \\
The set \,$\mathfrak{F}^1_{^\mathcal{V}}(r,\lambda)$ is relatively compact in $\mathfrak{F}^0(r,\lambda)$
in the following sense: \\[3mm]
\noindent Let $f^i:M^i\rightarrow \R^{m+1}$ be a sequence in $\mathfrak{F}^1_{^\mathcal{V}}(r,\lambda)$.
Then, after passing to a subsequence, there exist
an $f:M\rightarrow \R^{m+1}$ in $\mathfrak{F}^0(r,\lambda)$ and
a sequence of diffeomorphisms
$\phi^i:M\rightarrow M^i$, such that $f^i \circ \phi^i$ is
uniformly Lipschitz bounded and converges uniformly to $f$.
\end{theorem}
\vspace{5mm}
Here the Lipschitz bound for $f^i \circ \phi^i$ is shown with respect to the local representations of some
finite atlas of $M$. For these representations, we obtain a Lipschitz constant $L$ depending only on $\lambda$.
As an immediate consequence of Theorem \ref{compactness2} we deduce the following corollary: \\[-2mm]
\begin{cor} \label{corollaryfinit} There are only finitely many
manifolds in $\mathfrak{F}^1_{^\mathcal{V}}(r,\lambda)$ up to diffeomorphism.
\end{cor} \vspace{4mm}
The situation is slightly different when considering
$(r,\lambda)$-immersions in arbitrary codimension.
For the construction of the diffeomorphisms $\phi^i$ one uses a kind of projection in an
averaged normal direction $\nu$. In higher codimension,
the averaged normal $\nu$ cannot be constructed as in the case of
hypersurfaces. We will give an alternative construction involving a
Riemannian center of mass. However, for doing so we have to assume here
that $\lambda$ is not too large. Let $\mathfrak{F}^1_{^\mathcal{V}}(r,\lambda)$ and
$\mathfrak{F}^0(r,\lambda)$ be defined as above, but this time for
functions with values in $\R^{m+k}$ for a fixed \nolinebreak$k$.
We obtain the following theorem: \\
\begin{theorem} \label{compactness3} (Compactness of $(r,\lambda)$-immersions in arbitrary codimension) \\ \\[-4mm]
Let $\lambda\leq\frac{1}{4}$ .
Then $\mathfrak{F}^1_{^\mathcal{V}}(r,\lambda)$ is relatively compact in $\mathfrak{F}^0(r,\lambda)$ in
the sense of Theorem \ref{compactness2}.
\end{theorem}
\vspace{4mm}
As in Corollary \ref{corollaryfinit}, we deduce for $\lambda\leq\frac{1}{4}$ that
there are only finitely many manifolds
in $\mathfrak{F}^1_{^\mathcal{V}}(r,\lambda)$ up to diffeomorphism.
Surely, the bound $\lambda\leq\frac{1}{4}$ is not optimal; at the end of Section 6 we will discuss some possibilities how
to prove the theorem for bigger Lipschitz constant. \\ \\
In \cite{langer} and \cite{breuning1} any
sequence of immersions with $L^p$-bounded second fundamental form, $p>m$, is shown
to be also a sequence of $(r,\lambda)$-immersions (for some fixed $r$ and $\lambda$).
The same conclusion holds in many other situations, where the geometric data (such as
curvature bounds) ensure uniform graph representations with control over the
slope of the graphs. Hence it seems natural to unearth the compactness of
$(r,\lambda)$-immersions as a theorem on its own. In any general situation, where compactness
of immersions is desired (e.g.\ when considering convergence of geometric flows), only the
condition of Definition \ref{defgeneralizedimmersion} in Section 2 has to be verified. If in addition
some bound for higher derivatives of the graph functions is known (or for instance a $C^{0,\alpha}$-bound
for $Du$), with methods as in \cite{breuning1} one easily derives additional properties
of the limit, such as higher order differentiability or curvature bounds. Hence, Theorems
\ref{compactness2} and \ref{compactness3} can be seen as the most general kind of compactness
theorem in this context.
\\ \\ \\
\emph{Acknowledgement:} I would like to thank my advisor Ernst Kuwert for his support. Moreover
I would like to thank Manuel Breuning for proofreading
my dissertation \cite{breuningdiss}, where the results of this paper were established first.
\end{section}
\vspace{4mm}
\begin{section}{Definitions and preliminaries}
We begin with some general notations: For $n=m+k$ let $G_{n,m}$ denote the Grassmannian
of (non-oriented) $m$-dimensional subspaces of $\R^n$. Unless stated otherwise let
$B_{\!\varrho}$ denote the open ball in $\R^m$
of radius $\varrho>0$ centered at the origin. \\ \\
Now let $M$ be an $m$-dimensional manifold without boundary and $f:M\rightarrow\R^n$
a $C^1$-immersion. Let $q\in M$ and let $T_qM$ be the tangent space at $q$. Identifying vectors
$X\in T_qM$ with $f_\ast X\in T_{f(q)}\R^n$, we may consider $T_qM$ as an $m$-dimensional
subspace of $\R^n$. Let $(T_qM)^\bot$ denote the orthogonal complement of
$T_qM$ in $\R^n$, that is \beqnn
\R^n=T_qM\oplus(T_qM)^\bot \eeqnn
and $(T_qM)^\bot$ is perpendicular to $T_qM$. In this manner
we may define a tangent and a normal map
\\ \parbox{12cm}{\beqnn \tau_f:M&\rightarrow& G_{n,m}, \hspace{3cm}\\
q&\mapsto& T_qM,  \\ \text{and} \hspace{6.3cm} \\
 \nu_f:M&\rightarrow& G_{n,k}, \\
q&\mapsto& (T_qM)^\bot. \eeqnn}   \hfill  \parbox{8mm}{\beqn \label{notiontangentspace} \eeqn
\vspace{0.2cm} \beqn \label{normalnotion} \eeqn}
\\ \\ \\\textbf{The notion of an \boldmath$(r,\lambda)$\unboldmath-immersion:} \\ \\
We call a mapping $A:\R^n\rightarrow \R^n$ a \emph{Euclidean isometry}, if
there is a rotation $R\in \mathbb{SO}(n)$ and a translation $T\in\R^n$, such that
$A(x)=Rx+T$ for all $x\in \R^n$. \\ \\
\label{defeuclisom}For a given point $q \in M$ let $A_q: \mathbb{R}^n\rightarrow
\mathbb{R}^n$ be a Euclidean isometry,
which maps the origin to
$f(q)$, and the subspace $\mathbb{R}^m\times\{0\}\subset
\mathbb{R}^m \times \mathbb{R}^k$ onto $f(q)+\tau_f(q)$.
Let $\pi:\R^n\rightarrow \R^m$ be the standard projection onto the
first $m$ coordinates. \\ \\
Finally let $U_{r,q}\subset M$ be the $q$-component
of the set \label{qcomponentof} $(\pi \circ A_q^{-1} \circ f)^{-1}(B_r)$.
Although the isometry $A_q$ is not uniquely determined, the set
$U_{r,q}$ does not depend on the choice of $A_q$.\\ \\
We come to the central definition (as first defined in \cite{langer}):
\begin{defi} \label{definition1rlambda} An immersion $f$ is called an
$(r,\lambda)$-immersion, if for each point $q\in M$ the set
  $A_q^{-1}\circ f(U_{r,q})$ is the graph of a differentiable function
  $u:B_r\rightarrow \R^k$ with $\|Du\|_{C^{0}(B_{r})} \leq \lambda$. \end{defi}
\vspace{6mm}
Here, for any $x\in B_r$ we have $Du(x)\in \R^{k\times m}$. In order to define the
$C^0$-norm for $Du$, we have to fix a matrix norm for $Du(x)$. Of course all norms on
$\R^{k\times m}$ are equivalent, therefore our results are true for any norm
(possibly up to multiplication by some positive constant). Let us agree upon \vspace{-4mm}
\beqnn \|A\|\:=\:\Biggl(\hspace{0.5mm}\sum_{j=1}^m|a_j|^2\Biggr)^{\!\frac{1}{2}} \eeqnn
for $A=(a_1,\ldots,a_m)\in \R^{k\times m}$. For this norm we have
$\|A\|_{\text{op}}\:\leq\: \|A\|$ for any $A\in \R^{k\times m}$
and the operator norm $\|\cdot\|_{\text{op}}$.
Hence the bound $\|Du\|_{C^0(B_r)}\leq \lambda$ directly implies that $u$
is $\lambda$-Lipschitz.
Moreover the norm $\|Du\|_{C^{0}(B_{r})}$ does not depend on the choice of the
isometry \nolinebreak$A_q$. \pagebreak \\
\textbf{The notion of a generalized \boldmath$(r,\lambda)$\unboldmath-immersion:} \\ \\
For any $(r,\lambda)$-immersion $f:M\rightarrow \R^n$ and any $q\in M$,
we have a local graph representation over the affine tangent space
$f(q)+\tau_f(q)$. It is natural to extend this definition to immersions with local
graph representations over other appropriately chosen $m$-spaces in $\R^n$. \\ \\
For a given $q \in M$ and a given $m$-space $E\in G_{n,m}$ let $A_{q,E}: \mathbb{R}^n\rightarrow
\mathbb{R}^n$ be a Euclidean isometry, which maps the origin to
$f(q)$, and the subspace $\mathbb{R}^m\times\{0\}\subset
\mathbb{R}^m \times \mathbb{R}^k$ onto $f(q)+E$. \\ \\
Let $U_{r,q}^E\subset M$ be the $q$-component
of the set $(\pi \circ A_{q,E}^{-1} \circ f)^{-1}(B_r)$.
Again the isometry $A_{q,E}$ is not uniquely determined but the set
$U_{r,q}^E$ does not depend on the choice of $A_{q,E}$.
\\[-2mm]
\begin{defi} \label{defgeneralizedimmersion}
An immersion $f$ is called a generalized
$(r,\lambda)$-immersion, if for each point \linebreak$q\in M$ there is an $E=E(q)\in G_{n,m}$, such that the set
 $A_{q,E}^{-1}\circ f(U_{r,q}^E)$ is the graph of a differentiable function
 $u:B_r\rightarrow \R^k$ with $\|Du\|_{C^{0}(B_{r})} \leq
\lambda$.
\end{defi}
\vspace{1mm}
Obviously every $(r,\lambda)$-immersion is a generalized $(r,\lambda)$-immersion, as we can choose
$E(q)=\tau_f(q)$ for any $q\in M$. \\ \\
For fixed dimension $m$ and codimension $k$ we denote by $\mathfrak{F}^1(r,\lambda)$ the set of
generalized $(r,\lambda)$-immersions $f:M \rightarrow \R^{m+k}$ with $0\in f(M)$,
where $M$ is any compact $m$-manifold
without boundary. For $\mathcal{V}>0$ we denote by $\mathfrak{F}_{^\mathcal{V}}^1(r,\lambda)$ the
set of all immersions in $\mathfrak{F}^1(r,\lambda)$ with $\vol(M)\leq \mathcal{V}$.
Here the volume of $M$ is measured with respect to the volume measure induced by the metric
$f^\ast g_{\eucl}$.
Note that $M$ is \emph{not} fixed in these sets (in order to obtain a set in a strict
set theoretical sense one may consider every manifold as embedded in $\R^N$ for an  $N=N(m)$).
The condition
$0\in f(M)$ can be weakened in many applications to $f(M)\cap K\neq \emptyset$
for a compact set $K\subset \R^{m+k}$.
\\ \\
The notion of a generalized $(r,\lambda)$-immersion has one major advantage: As the definition does
not make use of the existence of a tangent space, it allows us to define similar notions
for functions into $\R^n$ which are not immersed. For a given $E\in G_{n,m}$ the set
$U_{r,q}^E$ can be defined for \emph{any} continuous function $f:M\rightarrow \R^n$. Moreover the
condition $\|Du\|_{C^{0}(B_{r})} \leq \lambda$ in the smooth case corresponds
to a Lipschitz bound of the function $u$. Hence the following
definition can be seen as the natural generalization to continuous functions: \\[-1mm]
\begin{defi} \label{rlfunction}
A continuous function $f$ is called an
$(r,\lambda)$-function, if for each point $q\in M$ there is an $E=E(q)\in G_{n,m}$, such that the set
 $A_{q,E}^{-1}\circ f(U_{r,q}^E)$ is the graph of a Lipschitz continuous function
 $u:B_r\rightarrow \R^k$ with Lipschitz constant $\lambda$.
\end{defi}
\vspace{3mm}
We additionally assume here, that $E$ can be chosen such that $f$ is injective on $U_{r,q}^E$.
This property is \emph{not} implied by the preceding definition, if one reads the latter word for word. \\ \\
We shall always consider $(r,\lambda)$-functions defined on compact topological
manifolds (without boundary).  Using the local Lipschitz graph representation,
any such manifold can be endowed with an atlas with
bi-Lipschitz change of coordinates. If the Lipschitz constant of the graphs is sufficiently
small (and hence the coordinate changes are almost isometric with
bi-Lipschitz constant close to $1$), by the results in \cite{karcher} there exists even a smooth atlas.
In our case, the limit manifold both in Theorem \ref{compactness2} and \ref{compactness3}
will be smooth.
\\ \\
Finally, we define the set $\mathfrak{F}^0(r,\lambda)$ by replacing generalized $(r,\lambda)$-immersions
in $\mathfrak{F}^1(r,\lambda)$ by $(r,\lambda)$-functions. \pagebreak \\
\textbf{Geometry of Grassmann manifolds} \\ \\
For $k,n\in \N$ with $0<k<n$ let $G_{n,k}$ again be the set of (non-oriented)
$k$-dimensional subspaces of $\R^n$. \\ \\
The set $G_{n,k}$ may be endowed with the structure of a differentiable $k(n-k)$-dimensional manifold,
see e.g.\ \cite{leed}. Moreover there is a Riemannian metric $g$ on $G_{n,k}$ being
invariant under the action of $\mathbb{O}(n)$ in $\R^n$. It is unique up to multiplication by a positive
constant (and --- again up to multiplication by a positive constant --- the only metric being invariant
under the action of $\mathbb{SO}(n)$ in $\R^n$ except for the case $G_{4,2}$). For more details we refer
the reader to \cite{leichtweiss}. \\ \\
In general, if $(M,g)$ is a Riemannian manifold, the induced distance on $M$ is defined by \\
\parbox{12cm}{\beqnn
d(p,q)=\inf \{L(\gamma)\:| &\gamma&\hspace{-3mm}:[a,b]\rightarrow M \text{ piecewise smooth curve with }
\\ &\gamma&\hspace{-3.5mm}(a)=p, \:\gamma(b)=q\}. \eeqnn} \hfill
\parbox{8mm}{\beqn \label{induceddistance} \eeqn} \\
Here $L(\gamma):=\int_{a}^b \lvert\frac{d\gamma}{dt}(t)\rvert\,dt$ denotes the length of $\gamma$.
If $M$ is complete, by the Theorem of Hopf-Rinow any two points $p,q\in M$ can be
joined by a geodesic of length $d(p,q)$. This applies to the Grassmannian as $G_{n,k}$ is complete. \\ \\
Now suppose that $E, G\in G_{n,k}$ are two close $k$-planes;
this means that the projection of each onto the
other is non-degenerate.
Applying a transformation to principal axes, there are orthonormal bases
$\{v_1,\ldots,v_k\}$ of $E$ and $\{w_1,\ldots,w_k\}$ of $G$ such that \beqnn
\langle v_i,w_j\rangle= \delta_{ij}\hspace{0.3mm}\cos \theta_i \hspace{5mm} \text{with } \theta_i\in
\left[0,\frac{\pi}{2}\right) \eeqnn
for $1\leq i,j \leq k$. For given $k$-spaces $E$ and $G$,
the $\theta_1,\ldots,\theta_k$ are uniquely determined (up to the order)
and called the \emph{principal angles} between $E$ and $G$.
Under all metrics on
$G_{n,k}$ being invariant under the action of $\mathbb{O}(n)$,
there is exactly one metric $g$ with \beqnn
d(E,G)=\left(\:\sum_{i=1}^k \theta_i^2\:\right)^{\frac{1}{2}} \eeqnn
for all close $k$-planes $E$ and $G$, where $d$ denotes the distance
corresponding to $g$, and $\theta_1,\ldots,\theta_k$ the principal angles
between $E$ and $G$ as defined above;
see \cite{borisenko} and the references given there.
We shall always use this distinguished metric. \\ \\
We will need the following estimate for the sectional curvatures of a Grassmannian:
\begin{lemma} \label{sectcurgrass} Let $\max\{k,n-k\}\geq2$. Let $K(\cdot,\cdot)$
denote the sectional curvature of $G_{n,k}$ and let
$X,Y \in T_{\!_P}G_{n,k}$ be linearly independent tangent vectors for a $P\in G_{n,k}$. Then \beqnn
0\leq K(X,Y)\leq 2. \eeqnn \end{lemma}
\textbf{Proof:} \\
For $\min\{k,n-k\}=1$ all sectional curvatures are constant with $K(X,Y)=1$. For a proof see
\cite{leichtweiss}, p.\ 351. For $\min\{k,n-k\}\geq 2$ we have $0\leq K(X,Y)\leq 2$ by
\cite{wong}, Theorem 3. \hfill $\square$ \\ \\
The injectivity radius of $G_{n,k}$ is $\frac{\pi}{2}$ (see \cite{borisenko}, p.\ 53).
A subset $U$ of a Riemannian manifold $(M,g)$ is said to be \emph{convex,} if and only if
for each $p,q \in U$ the shortest geodesic from $p$ to $q$ is unique in $M$ and lies
entirely in $U$.
For the Grassmannian $G_{n,k}$,
any open Riemannian ball $B_\varrho(P)$ around $P\in G_{n,k}$ with
$\varrho< \frac{\pi}{4}$ is convex; see \cite{gromoll}, p.\ 228. \\ \\[2mm]
\textbf{The Riemannian center of mass} \\ \\
The well-known Euclidean center of mass may be generalized to
a \emph{Riemannian center of mass} on Riemannian manifolds. This was introduced by
K. Grove and H. Karcher in \cite{grove}. A simplified treatment is given in \cite{karcher}.
See also \cite{jost}. We like to give a short sketch of this concept.
\\ \\
Let $(M,g)$ be a complete Riemannian manifold with induced distance
$d$ as in (\ref{induceddistance}). Let $\mu$ be a probability measure on $M$,
i.e.\ a nonnegative measure with \beqnn
\mu(M)=\int_M d\mu=1. \eeqnn
Let $q$ be a point in $M$ and $B_\varrho=B_\varrho(q)$
a convex open ball of radius $\varrho$ around
$q$ in $M$. Suppose \beqnn
\text{spt } \mu \subset B_\varrho\text{\emph{,}} \eeqnn
where spt $\mu$ denotes the support of $\mu$.
We define a function \beqnn
P:\overline{B}_\varrho&\rightarrow& \R, \\
P(p)&=& \int_M d(p,x)^2\,d\mu(x). \eeqnn
\begin{defi} A \:$q\in \overline{B}_\varrho$ is called a center of mass for $\mu$ if \beqnn
P(q)=\inf_{p\in \overline{B}_\varrho} \int_M d(p,x)^2\,d\mu(x). \eeqnn \end{defi} \vspace{2mm}
The following theorem asserts the existence and uniqueness of a center of mass:
\begin{theorem} \label{centerofmass}
If the sectional curvatures of $M$ in $B_\varrho$ are at most $\kappa$ with $0<\kappa<\infty$
and if $\varrho$ is small enough such that $\varrho<\frac{1}{4}\pi \kappa^{-1/2}$, then $P$ is
a strictly convex function on $B_\varrho$ and has a unique
minimum point in $\overline{B}_\varrho$ which lies in
$B_\varrho$ and is the unique center of mass for $\mu$. \end{theorem}
\textbf{Proof:} \\
See \cite{karcher}, Theorem 1.2 and the following pages there. \hfill $\square$ \\ \\
\noindent In the preceding theorem, we do not require the bound $\kappa$ to
be attained; in particular all sectional curvatures are also allowed to be less than or equal to $0$.
The same applies to the following lemma: \begin{lemma} \label{distcenterofmass}
Assume that the sectional curvatures of $M$ in $B_\varrho$ are at most $\kappa$ with $0<\kappa<\infty$
and $\varrho<\frac{1}{4}\pi \kappa^{-1/2}$.
Let $\mu_1,\mu_2$ be two probability measures on $M$ with spt $\mu_1\subset B_\varrho$,
spt $\mu_2\subset B_\varrho$ with centers of mass $q_1,q_2$ respectively. Then for a
universal constant $C=C(\kappa,\varrho)<\infty$ \beqnn
d(q_1,q_2)\leq C\int_M d(q_2,x)\,d|\mu_1-\mu_2|(x), \eeqnn
where $|\mu_1-\mu_2|$ denotes the total variation measure of the signed measure $\mu_1-\mu_2$.
\end{lemma}
\textbf{Proof:} \\
Let $P_i(p)=\frac{1}{2}\int_M d(p,x)^2\,d\mu_i(x)$ for $i=1,2$. By Theorem 1.5.1 in \cite{karcher}, with
\beqn \label{Ckappar} C=C(\kappa,\varrho):=1+(\kappa^{1/2}\varrho)^{-1}\tan(2 \kappa^{\,1/2}\varrho), \eeqn
we have for all $y\in B_\varrho$ the estimate
\beqnn d(q_1,y)\leq C \:|\text{grad } P_1(y)|. \eeqnn
Using spt $\mu_i\subset B_\varrho$\,,\, by Theorem 1.2 in \cite{karcher} we have \beqn \label{gradexpformula}
\text{grad } P_i(y)=-\int_{B_\varrho} \exp_y^{-1}(x) \,d\mu_i(x), \eeqn
where $\exp_y^{-1}:B_\varrho\rightarrow T_yM$ is considered as a vector valued function. \\ \\
Moreover, as $q_2$ is a center of mass, \beqnn
\text{grad } P_2(q_2)=0. \eeqnn
Then with the argumentation of \cite{jost}, Lemma 4.8.7 (where manifolds
of nonpositive sectional curvature are considered) we have \beqnn
d(q_1,q_2) &\leq& C\:|\text{grad } P_1(q_2)| \\
&=& C\:\biggl\lvert\int_{B_\varrho} \exp_{q_2}^{-1}(x) \,d\mu_1(x)\biggr\rvert \\
&=& C\:\biggl\lvert\int_{B_\varrho}\exp_{q_2}^{-1}(x) \,d\mu_1(x)-
\int_{B_\varrho}\exp_{q_2}^{-1}(x) \,d\mu_2(x)\biggr\rvert \\
&\leq& C\int_M d(q_2,x)\,d|\mu_1-\mu_2|(x), \eeqnn
where we used $|\exp_{q_2}^{-1}(x)|=d(q_2,x)$ and
spt $\mu_i\subset B_\varrho$ in the last line. \hfill $\square$ \\ \\ \\
\textbf{Basics for the proof} \\ \\
We like to fix some further notation and to deduce some basic facts that are
needed in the proof. \\ \\
First of all let us simplify the notation. For a given $(r,\lambda)$-immersion $f:M\rightarrow \R^{m+1}$
and for every $q\in M$ we can choose an $E_q \in G_{m+1,m}$
with the properties of Definition \ref{defgeneralizedimmersion}.
This yields a mapping $\mathcal{E}:M\rightarrow G_{m+1,m}$, $q \mapsto E_q$. For every
$(r,\lambda)$-immersion we choose and fix such a mapping $\mathcal{E}$.
So every given $(r,\lambda)$-immersion $f$ can be thought of as a pair $(f,\mathcal{E})$, even
if $\mathcal{E}$ is not explicitly mentioned in the notation. With
$A_{q,\mathcal{E}(q)}$ and $U_{r,q}^{\mathcal{E}(q)}$ as in
Definition \ref{defgeneralizedimmersion},
we set
\beqnn
A_q:=A_{q,\mathcal{E}(q)} \eeqnn
and for $0<\varrho \leq r$
\beqnn
U_{\varrho,q}:=U_{\varrho,q}^{\mathcal{E}(q)}. \eeqnn
In fact this means that $A_q$ and $U_{\varrho,q}$ also depend on $\mathcal{E}(q)$.
However, all properties shown below for $U_{\varrho,q}$ are true for any
admissible choice of $\mathcal{E}$.
\vspace{6mm} \\
As an analogue to Lemma 3.1 in \cite{langer} we obtain the following statement,
where $f$ is assumed to be a \emph{generalized}
$(r,\lambda)$-immersion here: \vspace{2.5mm}
\begin{lemma}\label{intersect} Let $f:M\rightarrow \R^{m+1}$ be an
$(r,\lambda)$-immersion and $p,q\in M$. \begin{itemize}
\item[a)] If\, $0<\varrho \leq r$ and $p \in U_{\varrho,q}$, then
$|f(q)-f(p)|<(1+\lambda)\varrho$.
\item[b)] If\, $0<\varrho \leq r$ and
$\delta=[3(1+\lambda)]^{-1}\varrho$ and $U_{\delta,q}\cap U_{\delta,p}\neq \emptyset$,
then $U_{\delta,p}\subset U_{\varrho,q}$. \end{itemize}
\end{lemma} \pagebreak
\noindent \textbf{Proof:}
\begin{itemize}
\item[a)] Pass to the graph representation, use the bound on the
$C^0$-norm of the derivative of the graph and the triangular inequality.
\item[b)] Let $x \in U_{\delta,p}$ and $y \in U_{\delta,q}\cap
U_{\delta,p}$. With $\varphi_q:=\pi\circ A_{q}^{-1}\circ
f$ we have \beqnn
|\varphi_q(x)| &\leq& |f(x)-f(q)| \\
&\leq& |f(x)-f(p)|+|f(p)-f(y)|+|f(y)-f(q)| \\
&<& 3(1+ \lambda)\delta \\
&=& \varrho. \eeqnn Hence
$U_{\delta,p} \subset \varphi_{q}^{-1}(B_\varrho)$. But $U_{\delta,p}\cup
U_{\delta,q}$ is a connected set containing $q$, hence included in the
$q$-component of
$\varphi_q^{-1}(B_\varrho)$, that is in $U_{\varrho,q}$. We conclude
$U_{\delta,p}\subset U_{\varrho,q}$. \hfill $\square$
\end{itemize} \vspace{8mm}
\noindent Now let $r,\lambda >0$ be given. For $l \in \N_0$ define
$\delta_l:=[3(1+\lambda)]^{-l}r$.
For an $(r,\lambda)$-immersion $f:M\rightarrow \R^{m+1}$, by
Lemma \ref{intersect} b) we have the following important property:
\beqn \label{deltalinclusion1}
\text{If } \,p,q \in M \text{ and } \,U_{\delta_{l+1},q}\cap U_{\delta_{l+1},p}\neq \emptyset, \text{ then }\,
U_{\delta_{l+1},p}\subset U_{\delta_{l},q}. \eeqn \\
If $f:M\rightarrow \R^{m+1}$ is an $(r,\lambda)$-immersion and $p\in M$, we may use the local graph
representation to conclude that
the set $f(U_{r,p})$ is homeomorphic to the ball $B_r$.
Hence we may choose a continuous unit normal
$\nu_p:U_{r,p}\rightarrow \mathbb{S}^m$
with respect to $f|U_{r,p}$. If $q\in M$ is another point and $\nu_q:U_{r,q}\rightarrow \mathbb{S}^m$
a continuous unit normal on $U_{r,q}$, we note that $\nu_p$ and $\nu_q$
do not necessarily
coincide on $U_{r,p}\cap U_{r,q}$. However, we have the following statement: \\
\begin{lemma} \label{normalsignintersection}
Let $f:M\rightarrow \R^{m+1}$ be an $(r,\lambda)$-immersion and $p,q\in M$. Let
$\nu_p:U_{\delta_1,p}\rightarrow \mathbb{S}^m$, $\nu_q:U_{\delta_1,q}\rightarrow \mathbb{S}^m$ be
continuous unit normals.
Suppose $U_{\delta_1,p}\cap U_{\delta_1,q}\neq \emptyset$.
Then exactly one of the following two statements is true: \begin{itemize}
\item $\nu_p(x)=\nu_q(x)$ \hspace{2.8mm}for every $x\in U_{\delta_1,p}\cap U_{\delta_1,q}$,
\item $\nu_p(x)=-\nu_q(x)$ for every $x\in U_{\delta_1,p}\cap U_{\delta_1,q}$.
\end{itemize}
\end{lemma} \vspace{1.5mm}
\textbf{Proof:} \\
Choose a $\xi\in U_{\delta_1,p}\cap U_{\delta_1,q}$. First suppose that
$\nu_p(\xi)=\nu_q(\xi)$. As $U_{r,p}$ is homeomorphic to $B_r$ and connected,
there are exactly two continuous unit normals on $U_{r,p}$. Let $\nu$ be the
one with $\nu(\xi)=\nu_p(\xi)$. Let
$W=\{x\in U_{\delta_1,p}:\nu(x)=\nu_p(x)\}$. Then $W$ is a nonempty subset
of the connected set $U_{\delta_1,p}$. Moreover $W$ is easily seen to be open and
closed in $U_{\delta_1,p}$. Therefore $W=U_{\delta_1,p}$ and $\nu_p=\nu$ on $U_{\delta_1,p}$. As
$U_{\delta_1,q}\subset U_{r,p}$ by (\ref{deltalinclusion1}), the preceding argumentation can
also be applied to $\nu_q$. With $\nu(\xi)=\nu_p(\xi)=\nu_q(\xi)$
we conclude $\nu_q=\nu$ on $U_{\delta_1,q}$. Hence $\nu_p=\nu=\nu_q$ on
$U_{\delta_1,p}\cap U_{\delta_1,q}$, as in the claim above.
If $\nu_p(\xi)=-\nu_q(\xi)$, a similar
argumentation yields $\nu_p=-\nu_q$ on $U_{\delta_1,p}\cap U_{\delta_1,q}$.
\hfill $\square$ \\[1mm]
\begin{rem}
The statement of the preceding lemma might seem to be obvious at first sight.
However one can think of a M\"{o}bius strip covered by two open sets $U$ and $V$,
each of which is homeomorphic to $B_r$, such that
$U\cap V$ has exactly two components. If we choose continuous unit normals
$\nu_1$, $\nu_2$ on $U,V$ respectively, we have $\nu_1=\nu_2$ on
one of the components, and $\nu_1=-\nu_2$ on the other.
Such a behavior of the normals is excluded by Lemma \ref{normalsignintersection}, irrespective whether
$U_{\delta_1,p}\cap U_{\delta_1,q}$ is connected or not.
\end{rem} \vspace{6mm}
We need the notion of a $\delta$-net:\\[-3mm]
\begin{defi}
Let $Q=\{q_1,\ldots,q_s\}$ be a finite set of points in $M$
and let $0<\delta<r$. We say that $Q$ is a $\delta$-net for $f$,
if $M=\bigcup \limits_{j=1}^{s}U_{\delta,q_j}$. \end{defi}
\vspace{0.3cm} \noindent Note that every $\delta$-net is also a
$\delta'$-net if $0<\delta<\delta'<r$. \\ \\ \\[-2mm]
The following statement is a bit stronger than Lemma 3.2 in \cite{langer}.
It bounds the number of elements in a $\delta$-net by an argumentation
similar to that in the proof of Vitali's covering theorem. Simultaneously, similarly
to Besicovitch's covering theorem, it gives a bound
(which does not depend on the volume) how often any fixed point in $M$
is covered by the net. More precisely, we have the following lemma:  \\
\begin{lemma} \label{cover} For $l\in\N$,
every $(r,\lambda)$-immersion on a compact $m$-manifold $M$ admits a
$\delta_l$-net $Q$ with \beqnn
|Q|&\leq& \delta_{l+1}^{-m}
\vol(M), \\
|\{q\in Q:p\in U_{\delta_{2},q}\}|&\leq& [3(1+\lambda)]^{(l+1)m} \hspace{7mm} \text{ for every fixed }\, p\in M. \eeqnn
\end{lemma}
\noindent \textbf{Proof:} \\
Let $q_1\in M$ be an arbitrary point. Assume we have found points
$\{q_1,\ldots,q_\nu\}$ in $M$ with the property
$U_{\delta_{l+1},q_j}\cap U_{\delta_{l+1},q_k}=\emptyset$ for $j\neq k$. Suppose
$U_{\delta_{l},q_1}\cup \ldots \cup U_{\delta_{l},q_\nu}$ does not cover $M$.
Then choose a point $q_{\nu+1}$ from the
complement. Then $U_{\delta_{l+1},q_k}\cap U_{\delta_{l+1},q_{\nu+1}}=\emptyset$ for $k\leq \nu$, as otherwise
$U_{\delta_{l+1},q_{\nu+1}}\subset U_{\delta_{l},q_{k}}$ by (\ref{deltalinclusion1}). As \beqnn
\vol(M)&\geq&\sum_{j=1}^s \vol(U_{\delta_{l+1},q_j}) \\
&\geq& \sum_{j=1}^s \mathcal{L}^m(B_{\delta_{l+1}}) \\
&\geq& s\delta_{l+1}^m, \eeqnn
this procedure yields after at most $\delta_{l+1}^{-m}\vol(M)$ steps a cover. \\ \\
For the second relation let $p\in M$. Let $Q=\{q_1,\ldots,q_s\}$ be the net that we found above.
Moreover let $Z(p)= \{q\in Q:p\in U_{\delta_{2},q}\}$.
By Lemma \ref{intersect} b) we have \beqnn
\bigcup_{q\in Z(p)}U_{\delta_2,q}\subset U_{\delta_1,p}. \eeqnn
Hence we may estimate as above \\ \parbox{14.5cm}{\beqnn
\vol(U_{\delta_{1},p})&\geq&\sum_{q\in Z(p)}\vol(U_{\delta_{l+1},q}) \\
&\geq& |Z(p)|\delta_{l+1}^m\mathcal{L}^m(B_1). \eeqnn} \hfill
\parbox{8mm}{\beqn \label{estimatede2net1} \eeqn} \\
As the immersion is an $(r,\lambda)$-immersion, we have \beqn \label{estimatede2net2}
\vol(U_{\delta_{1},p})\leq (1+\lambda)^m\delta_1^m\mathcal{L}^m(B_1). \eeqn
Combining (\ref{estimatede2net1}) and (\ref{estimatede2net2}), we estimate \beqnn
|Z(p)|&\leq& (1+\lambda)^m\delta_1^m\delta_{l+1}^{-m} \\
&=& 3^{lm}(1+\lambda)^{(l+1)m}, \eeqnn
which implies the statement. \hfill $\square$
\vspace{10mm} \\
\noindent We would like to emphasize that the second estimate in the preceding lemma
does not depend on the volume $\vol(M)$. This will
be necessary in order to obtain estimates for Lipschitz constants and for angles between different spaces depending
only on $\lambda$ but not on $\vol(M)$. \\ \\
\begin{defi} Let $f:M\rightarrow \R^{m+1}$ be an $(r,\lambda)$-immersion. Let $l\in \N$ and let
$Q=\{q_1,\ldots,q_s\}$ be a $\delta_l$-net for $f$. For $\iota\in \{0,1,\ldots,l\}$ and
$j \in \{1,\ldots,s\}$ we define \beqnn
Z_\iota(j):=\{1\leq k \leq s: U_{\delta_{\iota},q_{j}}\cap U_{\delta_{\iota},q_k}\neq \emptyset\}. \eeqnn
\end{defi}
\vspace{5mm}
\noindent For $\nu_1,\nu_2\in \R^{m+1}\backslash \{0\}$ let
$\sphericalangle(\nu_1,\nu_2)$ denote the non-oriented angle between $\nu_1$ and $\nu_2$, that is
\beqnn 0\: \leq\: &\sphericalangle(\nu_1,\nu_2)& \:\leq\: \pi, \\
&\sphericalangle(\nu_1,\nu_2)&=\arccos \frac{\langle \nu_1,\nu_2 \rangle}{|\nu_1||\nu_2|}.
\eeqnn \\
We consider the metric space $(\mathbb{S}^m,d)$, where
$\mathbb{S}^m\subset \R^{m+1}$ is the $m$-dimensional unit sphere and $d$ the intrinsic metric on
$\mathbb{S}^m$, that is
\beqn \label{metricsn}
d(\cdot,\cdot)=\,\sphericalangle(\cdot,\cdot). \eeqn
For $A\subset \mathbb{S}^m$ and $x\in \mathbb{S}^m$ let
dist$(x,A)=\inf\{d(x,y):y \in A\}$. For $\varrho>0$ let
$B_\varrho(A)=\linebreak\{x \in \mathbb{S}^m:\text{dist}(x,A)<\varrho\}$.
Moreover let $\mathcal{S}\subset \mathcal{P}(\mathbb{S}^m)$ denote the set of closed
nonempty subsets of $\mathbb{S}^m$. We denote by $d_\mathcal{H}$ the Hausdorff metric on $\mathcal{S}$, given by \beqnn
d_\mathcal{H}:\;\mathcal{S}\times \mathcal{S}\;&\rightarrow& \R_{\geq 0}, \\
(S_1,S_2)\!&\mapsto& \inf \{\varrho>0:\,S_1\subset B_\varrho(S_2),\, S_2\subset B_\varrho(S_1)\}. \eeqnn \\
We will need the following well-known version of the theorem of Arzel\`{a}-Ascoli for
the Hausdorff metric (see \cite{alt}, p.\ 125): \begin{lemma} \label{ArzelaHausdorff}
Let $(X,d)$ be a compact metric space and $\mathcal{A}$ the set of closed nonempty subsets of $X$. Then
$(\mathcal{A},d_\mathcal{H})$ is compact, i.e.\ every sequence in $\mathcal{A}$ has a subsequence that
converges to an element in $\mathcal{A}$. \end{lemma}
\vspace{5.5mm}
\noindent We \label{notationfortub} will have to estimate the size of some tubular neighborhoods.
To do this we need to introduce some more notation. Suppose we are given $\varrho>0$
and $u\in C^1(\overline{B}_\varrho)$
with $\|Du\|_{C^{0}(B_{\varrho})}\leq \lambda$. Moreover let
$T\in C^1(\overline{B}_\varrho,\R^{m+1})$ with $|T(x)|=1$ for all $x\in B_\varrho$.
Suppose that $T$ is $L$-Lipschitz for an $L$ with $0<L<\infty$. Let
$\omega:\overline{B}_\varrho\rightarrow G_{m+1,1}$, $q\mapsto \text{span }\{T(q)\}$.
Finally, let $\nu:B_\varrho\rightarrow \mathbb{S}^m$ be a continuous unit normal with respect
to the graph $x\mapsto(x,u(x))$.
We consider a vector bundle $E$ over $B_\varrho$, given by \beqnn
E=\{(x,y)\in B_\varrho \times \R^{m+1}: y \in \omega(x)\}. \eeqnn
For $\varepsilon>0$ let \beqnn
E^\varepsilon=\{(x,y)\in E: |y|<\varepsilon\}\subset E. \eeqnn
Moreover we define a mapping \\
\parbox{12cm}{\beqnn F:\;\;\;\;\; \,E\;\;&\rightarrow& \mathbb{R}^{m+1},\\
(x,y)&\mapsto&(x,u(x))+y, \eeqnn} \hfill
\parbox{8mm}{\beqn \eeqn} \\
where $y \in \omega(x)$. \\ \\

\begin{lemma}[Size of tubular neighborhoods] \label{tubularsigma} $ $ \\
 Let $\gamma<\frac{\pi}{2}$. With the notation as above,
assume that \beqn \label{angelTnu}
\sphericalangle(T(p),\nu(q)) \leq \gamma \hspace{5mm} \text{ for every }
p,q \in B_\varrho. \eeqn \vspace{-4mm} \\
Then the following is true: \begin{itemize}
\item[a)] For $\varepsilon=\frac{1}{L}\cos \gamma$ the mapping $F|E^\varepsilon$ is a
diffeomorphism onto an open neighborhood of $\{(x,u(x))\in \R^{m}\times \R:x\in B_\varrho\}$.
\item[b)] Let
$\sigma:=\min\!\left\{ \frac{\varrho}{2}\cos \gamma, \,\frac{\cos^2 \gamma}{2L(1+\lambda)}\!\right\}$.
Then
\beqnn
B_\sigma\left(\left\{(x,u(x))\in \R^m\times \R:x\in \overline{B}_{\frac{\varrho}{2}}\right\}\right)\subset
F(E^\varepsilon), \eeqnn
where $\varepsilon=\frac{1}{L}\cos \gamma$ as in part a) and
$B_\sigma(A)=\{x \in \R^{m+1}:\text{\emph{dist}}(x,A)<\sigma\}$ for $A\subset \R^{m+1}$ with
\emph{dist} the Euclidean
distance.
\end{itemize}
\end{lemma}
\vspace{11mm}
The trivial but long proof is carried out in detail in the appendix.
\\ \\ \\[9mm]
$\text{ }$ \hspace{-2cm} \setlength{\unitlength}{1.1cm}
\hspace{3cm}\begin{picture}(10,10)
\put(2,4){\line(1,0){8}} \put(5.95,3.925){$|$}
\put(9.92,3.925){)} \put(7.95,3.925){)}
\put(1.95,3.925){(} \put(3.95,3.925){(}
\put(10,3.65){$B_\varrho$}
\put(8,3.65){$B_{\frac{\varrho}{2}}$}

\linethickness{0.29mm}
\qbezier(2,9)(3,5)(6,4)
\qbezier(6,4)(9,3)(10,0.5)
\thinlines

\put(3.95,5.17){\tiny{(}}
\put(7.95,3.0){\tiny{)}}

\linethickness{0.29mm}
\qbezier(4.35,5.57)(4.9,4.9)(6.35,4.325)
\qbezier(6.35,4.325)(7.7,3.8)(8.35,3.4)
\qbezier(3.65,4.87)(4.2,4.2)(5.65,3.625)
\qbezier(5.65,3.625)(7,3.1)(7.65,2.7)

\qbezier(4.35,5.57)(3.3,5.77)(3.65,4.87)
\qbezier(8.35,3.4)(8.7,2.454)(7.65,2.7)
\thinlines

\put(2,9){\line(6,1){2.465975}}
\put(2,9){\line(-6,-1){2.465975}}
\put(2.3,8){\line(1,0){2.5}}
\put(2.3,8){\line(-1,0){2.5}}
\put(2.7,7){\line(6,-1){2.465975}}
\put(2.7,7){\line(-6,1){2.465975}}
\put(3.3,6){\line(6,1){2.465975}}
\put(3.3,6){\line(-6,-1){2.465975}}
\put(4.23,5){\line(1,0){2.5}}
\put(4.23,5){\line(-1,0){2.5}}
\put(6,4){\line(6,1){2.465975}}
\put(6,4){\line(-6,-1){2.465975}}
\put(8.15,2.9){\line(1,0){2.5}}
\put(8.15,2.9){\line(-1,0){2.5}}
\put(9.35,1.7){\line(5,-1){2.45145}}
\put(9.35,1.7){\line(-5,1){2.45145}}
\put(10,0.5){\line(6,1){2.465975}}
\put(10,0.5){\line(-6,-1){2.465975}}

\qbezier(3.95,9.33)(4.2,8.8)(4.3,8)
\qbezier(4.3,8)(4.4,7.1)(4.67,6.67)
\qbezier(4.67,6.67)(4.93,6.37)(5.26,6.33)
\qbezier(5.26,6.33)(5.63,6.3)(6.23,5)
\qbezier(6.23,5)(6.6,4.4)(7.97,4.34)
\qbezier(7.97,4.34)(9.2,4.1)(10.15,2.9)
\qbezier(10.15,2.9)(11.1,1.7)(11.3,1.3)
\qbezier(11.3,1.3)(11.64,0.83)(11.96,0.83)
\qbezier(0.025,8.66)(0.15,8.3)(0.3,8)
\qbezier(0.3,8)(0.65,7.6)(0.72,7.325)
\qbezier(0.72,7.325)(0.85,7)(0.9,6.7)
\qbezier(0.9,6.7)(1,6)(1.33,5.67)
\qbezier(1.33,5.67)(1.9,5.3)(2.23,5)
\qbezier(2.23,5)(3,4.4)(3.4,4)
\qbezier(3.4,4)(3.78,3.73)(4.023,3.683)
\qbezier(4.023,3.683)(5.245,3.46)(6.15,2.9)
\qbezier(6.15,2.9)(6.72,2.6)(7.38,2.1)
\qbezier(7.38,2.1)(7.9,1.8)(7.87,1.2)
\qbezier(7.87,1.2)(7.83,0.6)(8.02,0.165)

\put(5,9.14){\vector(-4,-1){2.75}}
\put(5.05,9.1){\small{$\{(x,u(x))\in \R^{m}\times \R:x\in B_\varrho\}$}}
\put(5.9,7.4){\vector(-1,0){2.2}}
\put(5.95,7.33){\small{$F(E^\varepsilon)$}}
\put(7.1,5.7){\vector(-2,-1){2}}
\put(7.15,5.63){\small{$B_\sigma\left(\left\{(x,u(x))\in \R^m\times \R:x\in \overline{B}_{\frac{\varrho}{2}}\right\}\right)$}}

\end{picture}
\vspace{-8mm}
\\ \begin{fig} \label{konzentration} Tubular neighborhood around a graph. \end{fig}
\vspace{1cm}
\noindent Finally we like to define again a metric for graph systems. First of all let
\beqnn \mathfrak{G}^s=\{(A_j,u_j)_{j=1}^s: \; A_j:
\mathbb{R}^{m+1}\rightarrow \mathbb{R}^{m+1} \text{ is a Euclidean isometry, } \;u_j \in
C^1(\overline{B}_r)\}. \eeqnn
Every Euclidean isometry $A:\R^{m+1}\rightarrow \R^{m+1}$ splits uniquely into a rotation
$R\in \mathbb{SO}(m+1)$ and a translation $T\in \R^{m+1}$. If $\|\cdot \|$ denotes
the operator norm and if $\Gamma=(A_j,u_j)_{j=1}^{s}\in \mathfrak{G}^s$,
$\tilde{\Gamma}=(\tilde{A}_j,\tilde{u}_j)_{j=1}^{s}\in
\mathfrak{G}^s$, we set
\label{metrik}
\\ \parbox{12cm}{ \beqnn \mathfrak{d}(\cdot,\cdot): \mathfrak{G}^s\times
\mathfrak{G}^s&\rightarrow& \mathbb{R}, \\
\mathfrak{d}(\Gamma,\tilde{\Gamma})&=&\sum
\limits_{j=1}^{s}(\|R_j-\tilde{R}_j\|+|T_j-\tilde{T}_j|+\|u_j-\tilde{u}_j\|_{C^{0}(B_{r})}).
\eeqnn} \hfill \parbox{8mm}{\beqn \label{metrikgraph} \eeqn} \\
This makes $(\mathfrak{G}^s,\mathfrak{d})$ a metric space.
\vspace{10mm}
\end{section}

\begin{section}{Transversality and tubular neighborhoods}
In this section we like to construct lines in $\R^{m+1}$,
that intersect each (appropriately restricted) immersion $f^i$ transversally --- even in the case, that the
Lipschitz constant $\lambda$ of the graph functions is large. This yields
local tubular neighborhoods around $f^i$ and is the crucial step in the proof. \\ \\
Let $r>0$ and $\lambda,\mathcal{V}<\infty$.
Let $f^i:M^i\rightarrow \R^{m+1}$ be a sequence of $(r,\lambda)$-immersions as in
Theorem \ref{compactness2}.
With Lemma \ref{cover} choose $\delta_5$-nets $Q^i=\{q_1^i,\ldots,q_{s^i}^i\}$
for $M^i$ with at most $\delta_6^{-m} \vol(M^i)$ elements
respectively and with \beqn \label{cardnet} \hspace{1.5cm}
|\{q\in Q^i:p\in U_{\delta_{2},q}^i\}|\leq [3(1+\lambda)]^{6m} \hspace{5mm}
\text{ for every fixed } p\in M^i. \eeqn
As $\vol(M^i)\leq \mathcal{V}$,
we may pass to a subsequence such that each net
has exactly $s$ points for a fixed $s \in \N$.\\ \\
For every $i\in \N$, $\iota\in\{0,1,\ldots,5\}$ and $j\in\{1,\ldots,s\}$ we have \beqnn
|Z_\iota^i(j)|\,\leq\, |\mathcal{P}(\{1,\ldots,s\})|\,=\,2^s. \eeqnn
Hence, by successively passing to subsequences, we may assume \beqn \label{Z4}
Z_\iota^i(j)=Z_\iota(j) \eeqn
for every $i,j$ and $\iota$ for fixed sets $Z_\iota(j)$. \\ \\
To simplify the notation, for $0<\varrho\leq r$ we set
$U_{\varrho,j}^i:=U_{\varrho,q_{j}^{i}}^{i}$. \\[2mm]
Moreover, we choose for every $i\in \N$ and every $j\in\{1,\ldots,s\}$ a continuous unit normal
$\nu_j^i:U_{r,j}^i\rightarrow \mathbb{S}^m$
with respect to $f^i|U_{r,j}^i$. Let these normal mappings be fixed from now on. \\ \\
For $S\subset \mathbb{S}^m$ let $\overline{S}$ be the closure of $S$ with respect
to the metric $d$ defined in (\ref{metricsn}).
We set \beqnn
S_j^i:=\overline{\nu_j^i(U_{\delta_{1},j}^i)}\subset \mathbb{S}^m. \eeqnn
For each
fixed $j$, this yields a sequence $(S_j^i)_{i\in \N}$ in $\mathcal{S}$. By Lemma \ref{ArzelaHausdorff}
we can pass successively to subsequences in order to obtain a sequence with
\beqnn
S_j^i \rightarrow S'_j\; \text{ in } (\mathcal{S},d_\mathcal{H}) \text{ as } i\rightarrow\infty \eeqnn
for each fixed $j\in \{1,\ldots,s\}$, where $S'_j\in\mathcal{S}$.
In particular for every $j$
\beqn \label{cauchy}
(S_j^i)_{i\in\N} \text{ is a Cauchy sequence in } (\mathcal{S},d_\mathcal{H}). \eeqn
By (\ref{cauchy}) we may choose another subsequence such that for every $j$
\beqn \label{cauchydist}
d_\mathcal{H}(S_j^k,S_j^l)\;<\; \frac{\pi}{4}-\frac{1}{2}\arctan \lambda \hspace{5mm}
\text{ for all } k,l \in \N. \eeqn
\vspace{3mm} \\
To each $q_j^i\in Q^i$ we may assign a neighborhood $U_{r,j}^i$, a Euclidean isometry $A_j^i$
and a differentiable function $u_j^i:B_r\rightarrow \R$ as in Definition \ref{defgeneralizedimmersion}.
This yields the corresponding graph systems
$\Gamma^i=(A_j^i,u_j^i)_{j=1}^s\in \mathfrak{G}^s$. As $\|Du_j^i\|_{C^0(B_r)}\leq \lambda$
and as $f^i(M^i)$ is uniformly bounded, a subsequence of
$(\Gamma^i)_{i\in\N}$ converges in $(\mathfrak{G}^s,\mathfrak{d})$.
In particular \beqn \label{cauchygraph}
(\Gamma^i)_{i\in\N} \text{ is a Cauchy sequence in } (\mathfrak{G}^s,\mathfrak{d}). \eeqn
Let constants $L,\gamma$ and $\sigma$ be defined by \beqn
 L&:=& [3(1+\lambda)]^{6m+4}r^{-1}, \label{def1L} \\[2.5mm]
\gamma&:=&\frac{\pi}{4}+\frac{1}{2}\arctan \lambda, \label{def2gamma} \\
\sigma&:=& \frac{\cos^2 \gamma}{2L(1+\lambda)}. \label{def3sigma}
\eeqn \\
By (\ref{cauchygraph}) we may pass to another subsequence such that
\beqn \label{dlesssigma}
\mathfrak{d}(\Gamma^k,\Gamma^l)\;<\;[3(1+\lambda)(1+r)]^{-1}\sigma \hspace{8mm}
\text{ for all } k,l\in\N. \eeqn \\
For $i=1$ we sometimes suppress the index $1$ and write for instance $q_j$ and $u_j$ instead of $q_j^1$
and $u_j^1$. For the immersion $f^1$, let $\mathcal{E}^1:M^1\rightarrow G_{m+1,m}$ be a mapping
as explained in the beginning of Chapter 3.1.
We set $E_j:=\mathcal{E}^1(q_j^1)\in G_{m+1,m}$ (this means $E_j$ is an $m$-space
for the point $q_j^1\in M^1$ as in Definition \ref{defgeneralizedimmersion}). \\ \\ \\
Our next task is to find a mapping $\omega:M^1\rightarrow G_{m+1,1}$, which defines the direction in which
we project from $f^1(M^1)$ onto $f^i(M^i)$ in order to construct diffeomorphisms $\phi^i:M^1\rightarrow M^i$.
First we would like to give a local construction. In Lemma \ref{omegaisgloballywell} we will show that $\omega$
is even globally well-defined.
The construction is similar to that in \cite{langer}, but more involved. \\ \\
We choose a $C^\infty$-function
$g:\R_{\geq 0}\rightarrow \R$ with the following properties: \begin{itemize}
\item $g(t)=1$ for $t<\frac{\delta_1}{r}$,
\item $0\leq g(t) \leq 1$ for $t\in[\frac{\delta_1}{r},1]$,
\item $g(t)=0$ for $t>1$,
\item $-2 \leq g'(t) \leq 0$ for all $t>0$. \end{itemize}
We note that $\frac{\delta_1}{r}=[3(1+\lambda)]^{-1}\leq\frac{1}{3}$,
hence such a function $g$ exists.
\vspace{5mm} \\ Let \beqnn
Z:M^1&\rightarrow& \mathcal{P}(\{1,\ldots,s\}), \\
q&\mapsto& \{1 \leq k \leq s: q\in U_{\delta_{2},k}^1 \}. \eeqnn  \\
By (\ref{cardnet}) we have \beqn \label{cardZ}
|Z(q)|\;\leq\; [3(1+\lambda)]^{6m} \hspace{5mm} \text{ for every } q\in M^1. \eeqn
For every $k\in\{1,\ldots,s\}$ we choose a unit vector $w_k$ that is
perpendicular to the subspace $E_k$ defined above. Let
these vectors $w_1,\ldots,w_s$ be fixed from now on. \\ \\
Now let $j\in \{1,\ldots,s\}$,
$q \in U_{\delta_{3},j}^1$ and $k\in Z(q)$. Lemma \ref{intersect} b) yields \beqnn
U_{\delta_{1},j}^1\subset U_{r,k}^1. \eeqnn
In particular $f^1(U_{\delta_{1},j}^1)$ is the graph of a
$\lambda$-Lipschitz function on a subset of $E_k$. This implies \beqn \label{eitherorcase}
\text{either }\hspace{4mm} \sphericalangle(w_k,\nu_j^1(q_j)) \leq \arctan \lambda \hspace{4mm} \text{ or } \hspace{4mm}
\sphericalangle(-w_k,\nu_j^1(q_j)) \leq \arctan \lambda. \eeqn
Set \beqn \label{sign}
\nu_k:= \left\{ \begin{array}{c}
                        \;\;\;\,w_k, \hspace{4mm} \text{if } \sphericalangle(w_k,\nu_j^1(q_j)) \leq \arctan \lambda, \\
                        -w_k, \hspace{10.5mm} \text{otherwise}. \hspace{20mm} \text{ }
                      \end{array} \right.
\eeqn
If we replace the point $q_j$ by any other point $p\in U_{\delta_1,j}^1$, the relation (\ref{eitherorcase}) will still be true.
As $\nu_j^1$ is continuous and $U_{\delta_{1},j}^1$ is connected, we easily conclude \beqn \label{angle}
\sphericalangle(\nu_k,\nu_j^1(p))\leq \arctan \lambda \hspace{5mm} \text{ for every } p\in U_{\delta_{1},j}^1, \eeqn
where $\nu_k$ is the fixed vector defined in (\ref{sign}).
We finally define a function \beqnn
S:U_{\delta_{3},j}^1&\rightarrow& \R^{m+1}, \\
q&\mapsto& \sum_{k \in Z(q)}g\left(\frac{|f^1(q)-f^1(q_k)|}{\delta_{2}}\right)\nu_k. \eeqnn \\
\begin{lemma} \label{sestimate} The following inequalities hold:
\begin{itemize}
\item[a)] $|S(q)|\geq (1+\lambda)^{-1}$ for every $q \in U_{\delta_{3},j}^1$.
\item[b)]
$\sphericalangle(S(q),\nu_j^i(p))\leq \frac{\pi}{4}+\frac{1}{2}\arctan \lambda$\:
for every $q \in U_{\delta_{3},j}^1$ and every $p\in U_{\delta_{1},j}^i$.
\end{itemize}
\end{lemma}
\noindent \textbf{Proof:} \vspace{-2mm} \begin{itemize}
\item[a)] Let $q\in U_{\delta_{3},j}^1$.
As $Q^1$ is a $\delta_4$-net for $f^1$, there is a
$k\in \{1,\ldots,s\}$ with $q\in U_{\delta_{4},k}^1$.
By Lemma \ref{intersect} a) we have
$|f^1(q)-f^1(q_k)|<\delta_3$, hence \beqnn
\frac{|f^1(q)-f^1(q_k)|}{\delta_2}\;<\; \frac{\delta_3}{\delta_2}\;=\;\frac{\delta_1}{r}. \eeqnn
By the definition of $g$ this yields \beqnn
g\left(\frac{|f^1(q)-f^1(q_k)|}{\delta_2}\right)=1. \eeqnn
Now let $l\in Z(q)$.
By (\ref{angle}) we have $\sphericalangle(\nu_l,\nu_j^1(q))\leq \arctan \lambda$.
Hence \beqnn
\langle\nu_l,\nu_j^1(q)\rangle &=& |\nu_l||\nu_j^1(q)|\cos(\sphericalangle(\nu_l,\nu_j^1(q))) \\
&\geq& \cos\,(\arctan \lambda)\\
&=& (1+\lambda^2)^{-\frac{1}{2}} \\
&\geq& (1+\lambda)^{-1}. \eeqnn
We note that $q\in U_{\delta_{4},k}^1$ in particular implies $k\in Z(q)$. Finally we estimate \beqnn
|S(q)| &\geq& \langle S(q),\nu_j^1(q)\rangle \\
&=& g\left( \frac{|f^1(q)-f^1(q_k)|}{\delta_2}\right)
\langle\nu_k,\nu_j^1(q)\rangle
+\sum_{l\in Z(q)\setminus\{k\}}g\left( \frac{|f^1(q)-f^1(q_l)|}{\delta_2}\right)
\langle \nu_l,\nu_j^1(q)\rangle \\
&\geq& \left(1 + \sum_{l\in Z(q)\setminus\{k\}}g\left( \frac{|f^1(q)-f^1(q_l)|}{\delta_2}\right)\right)
(1+\lambda)^{-1} \label{angleestimate}\\
&\geq& (1+\lambda)^{-1}. \eeqnn
\item[b)]
Let $q \in U_{\delta_{3},j}^1$ and $p\in U_{\delta_{1},j}^i$.
By (\ref{cauchydist}) there is a $p'\in U_{\delta_{1},j}^1$ with \beqn \label{eq1}
\sphericalangle(\nu_j^1(p'),\nu_j^i(p))\;\leq \;
\frac{\pi}{4}-\frac{1}{2}\arctan \lambda. \eeqn
By (\ref{angle}), every $\nu_k$ with $k\in Z(q)$ lies in the cone \beqnn
C=\{v\in \R^{m+1}\setminus\{0\}:\sphericalangle(v,\nu_j^1(p'))\leq \arctan \lambda\}. \eeqnn
By the definition of $S$, also the non-zero vector $S(q)$ lies in $C$, i.e. \beqn \label{angletriangle}
\sphericalangle(S(q),\nu_j^1(p'))\leq \arctan \lambda. \eeqn
Using the triangular inequality, we conclude with (\ref{eq1}) and (\ref{angletriangle}) that \beqnn
\sphericalangle(S(q),\nu_j^i(p))\leq \frac{\pi}{4}+\frac{1}{2}\arctan \lambda. \eeqnn
\vspace{-0.9cm} \\ $\text{ }$ \hfill $\square$
\end{itemize} \vspace{5mm}
By Lemma \ref{sestimate} a) the mapping $S$ does not vanish on $U_{\delta_{3},j}^1$.
We define $T$ by normalizing $S$, that is
\beqnn
T:U_{\delta_{3},j}^1&\rightarrow& \R^{m+1}, \\
q&\mapsto& \frac{S(q)}{|S(q)|}. \eeqnn \\
Identifying $U_{\delta_3,j}^1$ with $B_{\delta_3}$ by means of the diffeomorphism
$\pi\circ A_j^{-1}\circ f^1:U_{\delta_3,j}^1\rightarrow B_{\delta_3}$,
we may consider $T$ and $S$ as mappings defined on the ball $B_{\delta_3}$. We show, that $T$
considered as mapping on $B_{\delta_3}$ is Lipschitz
with respect to the Euclidean norm: \\[0mm]
\begin{lemma} \label{TisLipschitz} The mapping $T:B_{\delta_{3}}\rightarrow \R^{m+1}$ is $L$-Lipschitz with
$L= [3(1+\lambda)]^{6m+4}r^{-1}$.
\end{lemma}
\noindent \textbf{Proof:} \\
Let $x,y\in B_{\delta_{3}}$. Then there are unique $p,q\in U_{\delta_{3},j}^1$ with
$\pi\circ A_j^{-1}\circ f^1(p)=x$, \;$\pi\circ A_j^{-1}\circ f^1(q)=y$. \\ \\
Let $k\in Z(p)\setminus Z(q)$. Then
$p\in U_{\delta_3,j}^1 \cap U_{\delta_2,k}^1$. Lemma \ref{intersect} b) implies
$U_{\delta_3,j}^1\subset U_{\delta_1,k}^1$, so in particular $q\in U_{\delta_1,k}^1$.
Now assume $|f^1(q)-f^1(q_k)|<\delta_2$. With $\varphi_k=\pi\circ A_k^{-1}\circ f^1$ this implies
$\varphi_k(q)\in B_{\delta_2}$. Hence $q\in U_{\delta_1,k}^1\cap \varphi_k^{-1}(B_{\delta_2})
=U_{\delta_2,k}^1$. But this contradicts $k \notin Z(q)$. Therefore
$|f^1(q)-f^1(q_k)|\geq \delta_2$ and hence
$g\left( \frac{|f^1(q)-f^1(q_k)|}{\delta_2}\right)=0$ by the definition of $g$. \\ \\
The same argument shows $g\left( \frac{|f^1(p)-f^1(q_l)|}{\delta_2}\right)=0$
for all $l\in Z(q)\setminus Z(p)$. \\ \\[2mm]
Using the preceding considerations, $\|g'\|_{C^0(\R_{\geq0})}\leq 2$ and $|Z(p)|\leq [3(1+\lambda)]^{6m}$,
$|Z(q)|\leq [3(1+\lambda)]^{6m}$, we estimate as follows: \beqnn
|S(x)-S(y)|&=&\Biggl\lvert \sum_{k\in Z(p)} g\left( \frac{|f^1(p)-f^1(q_k)|}{\delta_2}\right)\nu_k
-\sum_{l\in Z(q)} g\left( \frac{|f^1(q)-f^1(q_l)|}{\delta_2}\right)\nu_l\Biggr\rvert \\
&=& \Biggl\lvert\sum_{k\in Z(p)\cup Z(q)}\left[ g\left( \frac{|f^1(p)-f^1(q_k)|}{\delta_2}\right)
-g \left( \frac{|f^1(q)-f^1(q_k)|}{\delta_2}\right)\right] \nu_k\Biggr\rvert \\
&\leq& \sum_{k\in Z(p)\cup Z(q)} \|g'\|_{C^0(\R_{\geq0})}\left\lvert
\frac{|f^1(p)-f^1(q_k)|}{\delta_2}-\frac{|f^1(q)-f^1(q_k)|}{\delta_2}\right\rvert \\
&\leq& \sum_{k\in Z(p)\cup Z(q)} \frac{2}{\delta_2}|f^1(p)-f^1(q)| \\
&\leq& 4[3(1+\lambda)]^{6m+2}r^{-1}|(x,u_j(x))-(y,u_j(y))| \\
&\leq& 4[3(1+\lambda)]^{6m+2}r^{-1}(1+\lambda)|x-y|.
\eeqnn \\
By Lemma \ref{sestimate} a) we have $|S(z)|\geq(1+\lambda)^{-1}$
for every $z\in U_{\delta_{3},j}^1$. Hence \beqnn
|T(x)-T(y)|=\left\lvert\frac{S(x)}{|S(x)|} - \frac{S(y)}{|S(y)|}\right\rvert&\leq&
4[3(1+\lambda)]^{6m+2}(1+\lambda)^2r^{-1}|x-y|\\
&\leq& [3(1+\lambda)]^{6m+4}r^{-1}|x-y|. \eeqnn \nopagebreak
\vspace{-11.2mm} \\
$\text{ }$ \hfill $\square$ \vspace{2mm}
\begin{rem} Of course, $T$ is also Lipschitz as a mapping on $U_{\delta_3,j}^1$ with respect
to the metric induced by $f^1$. The estimate of the Lipschitz constant gets even better
in this case. Moreover, we note that in the preceding lemma $L$ depends on $r$.
However, we will see
that the Lipschitz constant of $f^i\circ \phi^i$ does not depend on $r$ in the end.
\end{rem} \vspace{5mm}
We set \beqnn
\omega:U_{\delta_{3},j}^1 &\rightarrow& G_{m+1,1}, \\
q&\mapsto& \text{span} \{S(q)\}, \eeqnn
which is well-defined as $S(q)\neq0$ by Lemma \ref{sestimate} a). \\ \\ \\
We like to explain how $\omega$ locally forms a tubular neighborhood around $f^1$: \\ \\
For that we consider the mapping \beqnn
g_k:U_{\delta_{2},k}^1&\rightarrow& \R, \\
q&\mapsto& g\left( \frac{|f^1(q)-f^1(q_k)|}{\delta_2} \right). \eeqnn
As $g$ is smooth and $g(t)=0$ for $t\geq 1$, it is easily seen that $g_k$ can be extended
to a smooth function $\bar{g}_k:M^1\rightarrow \R$ by setting $\bar{g}_k=0$ outside
$U_{\delta_{2},k}^1$. This implies that $S:U_{\delta_3,j}^1\rightarrow \R^{m+1}$ is
differentiable, even if the sum in the definition of $S$ depends on $Z(q)$.
Hence also $T=\frac{S}{|S|}$ is differentiable. Moreover Lemma \ref{TisLipschitz} says that
$T$ is $L$-Lipschitz with
$L=[3(1+\lambda)]^{6m+4}r^{-1}$
and by Lemma \ref{sestimate} b) we have
\beqnn
\sphericalangle(T(p),\nu_j^1(q))\leq \frac{\pi}{4}+\frac{1}{2}\arctan \lambda
\hspace{5mm} \text{ for all } p,q \in U_{\delta_{3},j}^1. \eeqnn
Finally, after a rotation and a translation, $f(U_{\delta_3,j}^1)$ may be written as the
graph of a $C^1$-function $u_j^1:B_{\delta_3}\rightarrow\R$. Let us introduce some more notation: \\ \\
We consider a vector bundle $\hat{E}_j$ over $U_{\delta_3,j}^1$, given by \beqnn
\hat{E}_j=\{(x,y)\in U_{\delta_{3},j}^1 \times \R^{m+1}: y \in \omega(x)\} \eeqnn
with bundle projection $\hat{\pi}$.
We may identify the zero section of $\hat{E}_j$ with $U_{\delta_3,j}^1$.
For $\varepsilon>0$ let \beqnn
E_j^\varepsilon=\{(x,y)\in \hat{E}_j: |y|<\varepsilon\}\subset \hat{E}_j. \eeqnn
Finally we define a mapping \\
\parbox{12cm}{\beqnn F_j:\;\;\;\;\; \,\hat{E}_j\;\;&\rightarrow& \mathbb{R}^{m+1},\\
(x,y)&\mapsto&f(x)+y, \eeqnn} \hfill
\parbox{8mm}{\beqn \eeqn} \\
where $y \in \omega(x)$. \\
\begin{lemma} \label{tubularF} Let $\varepsilon=\frac{1}{L}\cos \gamma$, where
$L$ and $\gamma$ are as in (\ref{def1L}), (\ref{def2gamma}). Then the following is true: \begin{itemize}
\item $F_j|E_j^\varepsilon$ is a diffeomorphism onto an open neighborhood of $f^1(U_{\delta_{3},j}^1)$,
\item $F_j|U_{\delta_{3},j}^1=f^1|U_{\delta_{3},j}^1$,
\item for each fibre $\hat{E}_q=\hat{\pi}^{-1}(q)$ it holds $F_j(\hat{E}_q)=\omega(q)$.
\end{itemize}
Moreover for $\sigma=\frac{\cos^2 \gamma}{2L(1+\lambda)}$ we have the inclusion \beqnn
B_\sigma(f^1(U_{\delta_{4},j}^1))\subset F_j(E_j^\varepsilon). \eeqnn
\end{lemma}
\noindent \textbf{Proof:} \\
This is just a reformulation of Lemma \ref{tubularsigma}. Note that
\beqnn \frac{\cos \gamma}{L(1+\lambda)}
\;<\; [3(1+\lambda)]^{-3}r
\;=\; \delta_3, \eeqnn
hence $\sigma=\min\left\{\frac{\delta_3}{2}\cos \gamma,
\frac{\cos^2 \gamma}{2L(1+\lambda)}\right\}
=\frac{\cos^2 \gamma}{2L(1+\lambda)}$.
 \hfill $\square$
\vspace{15mm} \\
\noindent Up to this point we have constructed for each $j\in\{1,\ldots,s\}$ a tubular neighborhood locally
around $f(U_{\delta_{3},j}^1)$. Since the mapping $S$ depends on $j$, we
should write more accurately $S_j$ instead of $S$. In the same way we should write
$\omega_j$ instead of $\omega$. However, we can show that $\omega$ is globally
well-defined. More precisely we have the following lemma: \\
\begin{lemma} \label{omegaisgloballywell} Let $j,k\in \{1,\ldots,s\}$. Then \beqnn
\omega_j=
\omega_k \hspace{4mm} \text{ on } \; U_{\delta_{3},j}^1\cap U_{\delta_{3},k}^1.
\eeqnn
In particular there is a smooth mapping $\omega:M^1\rightarrow G_{m+1,1}$ with
$\omega|U_{\delta_{3},j}^1=\omega_j$ for each $j\in\{1,\ldots,s\}$.
\end{lemma}
\noindent \textbf{Proof:} \\
Let $j,k\in\{1,\ldots,s\}$. For $q \in U_{\delta_{3},j}^1\cap U_{\delta_{3},k}^1$ we show
that either $S_j(q)=S_k(q)$ or $S_j(q)=-S_k(q)$, which implies the statement. \\ \\
Let $q \in U_{\delta_{3},j}^1\cap U_{\delta_{3},k}^1$ and $l\in Z(q)$.
Lemma \ref{intersect} b)
implies \beqnn
U_{\delta_{1},j}^1&\subset& U_{r,l}^1, \\
U_{\delta_{1},k}^1&\subset& U_{r,l}^1. \eeqnn
As in (\ref{eitherorcase}) we conclude
\beqnn
\left(\,\text{either } \hspace{4mm}\sphericalangle(w_l,\nu_j^1(q_j)) \leq \arctan \lambda \hspace{4mm} \text{ or } \hspace{4mm}
\sphericalangle(-w_l,\nu_j^1(q_j)) \leq \arctan \lambda\:\right) \eeqnn
and  \beqnn
\hspace{2mm}\left(\,\text{either } \hspace{3.75mm}\sphericalangle(w_l,\nu_k^1(q_k)) \leq \arctan \lambda \hspace{4mm} \text{ or } \hspace{4mm}
\sphericalangle(-w_l,\nu_k^1(q_k)) \leq \arctan \lambda\:\right). \eeqnn
We define vectors as in (\ref{sign}), the first time depending on $j$, the second time on $k$:
\beqnn
\nu_{j,l}:= \left\{ \begin{array}{c}
                        \;\;\;\,w_l, \hspace{4mm} \text{if } \sphericalangle(w_l,\nu_j^1(q_j)) \leq \arctan \lambda, \\
                        -w_l, \hspace{10mm} \text{otherwise}. \hspace{20mm} \text{ }
                      \end{array} \right., \\ \\
\nu_{k,l}:= \left\{ \begin{array}{c}
                        \;\;\;\,w_l, \hspace{4mm} \text{if } \sphericalangle(w_l,\nu_k^1(q_k)) \leq \arctan \lambda, \\
                        -w_l, \hspace{10mm} \text{otherwise}. \hspace{20mm} \text{ }
                      \end{array} \right..
\eeqnn
Then \beqnn
S_j(q)=\sum_{l\in Z(q)}g\left(\frac{|f^1(q)-f^1(q_l)|}{\delta_{2}}\right)\nu_{j,l}
\eeqnn and \beqnn
S_k(q)=\sum_{l\in Z(q)}g\left(\frac{|f^1(q)-f^1(q_l)|}{\delta_{2}}\right)\nu_{k,l}. \eeqnn
By Lemma \ref{normalsignintersection}, we have $\nu_j^1= \nu_k^1$ on $U_{\delta_1,j}^1\cap U_{\delta_1,k}^1$,
or $\nu_j^1= -\nu_k^1$ on $U_{\delta_1,j}^1\cap U_{\delta_1,k}^1$.
Let us first assume \beqn \label{supppose}
\nu_j^1=\nu_k^1 \:\text{ on } U_{\delta_1,j}^1\cap U_{\delta_1,k}^1. \eeqn \\ Since
$q\in U_{\delta_3,j}^1\cap U_{\delta_3,k}^1$, we conclude with Lemma \ref{intersect} b) \beqnn
U_{\delta_3,j}^1\subset U_{\delta_2,k}^1, \hspace{10mm} U_{\delta_3,k}^1\subset U_{\delta_2,j}^1, \eeqnn
in particular \beqn \label{subset1}
\{ q_j,q_k\}\;\subset\; U_{\delta_3,j}^1\cup U_{\delta_3,k}^1\;\subset\;
U_{\delta_1,j}^1\cap U_{\delta_1,k}^1. \eeqn
By (\ref{angle}) together with (\ref{subset1}) we have \beqn \label{ineqq1}
\sphericalangle(\nu_{j,l},\nu_j^1(q_k)) \leq \arctan \lambda, \eeqn
by (\ref{supppose}), (\ref{subset1}) and (\ref{ineqq1}) moreover \beqn \label{ineqq2}
\sphericalangle(\nu_{j,l},\nu_k^1(q_k)) \leq \arctan \lambda. \eeqn
We already know that $\nu_{j,l}=\nu_{k,l}$ or $\nu_{j,l}=-\nu_{k,l}$, thus (\ref{ineqq2})
allows us to conclude that \beqnn
\nu_{j,l}=\nu_{k,l}. \eeqnn
Since this is true for all $l\in Z(q)$, we conclude $S_j(q)=S_k(q)$ and hence
$\omega_j(q)=\omega_k(q)$. \\ \\
If $\nu_j^1=-\nu_k^1 \text{ on } U_{\delta_1,j}^1\cap U_{\delta_1,k}^1$, one similarly concludes
$\nu_{j,l}=-\nu_{k,l}$ for all $l\in Z(q)$. This implies $S_j(q)=-S_k(q)$ and hence again
$\omega_j(q)=\omega_k(q)$. \hfill $\square$
\end{section}
\begin{section}[Intersection points and definition of $\phi^i$]{Intersection points and definition of \boldmath$\phi^i$\unboldmath}
In this section we like to show that for $p\in M^1$ the line $f^1(p)+\omega(p)$
intersects each appropriately restricted immersion $f^i(M^i)$ in exactly one
point. Using this, we are able to give a definition of the mappings
$\phi^i:M^1\rightarrow M^i$. Each $\phi^i$ will be shown to be a diffeomorphism.
Moreover, it will be shown that $f^i\circ \phi^i$ is uniformly Lipschitz bounded. \\
\begin{lemma} \label{intersectionpoint} For $p\in U_{\delta_{3},j}^1$ the line $f^1(p)+\omega(p)$ intersects the
set $f^i(U_{\delta_{1},j}^i)$ in exactly one point. This point lies in $f^i(U_{\delta_{2},j}^i)$. \end{lemma}
\noindent \textbf{Proof:} \\
Let $p\in U_{\delta_3,j}^1$.
First we show that $f^1(p)+\omega(p)$ intersects $f^i(U_{\delta_2,j}^i)$.
By Lemma \ref{sestimate} b) we have
\beqn \label{anglepfixed}
\sphericalangle(T(p),\nu_j^i(q))\leq \frac{\pi}{4}+\frac{1}{2}\arctan \lambda \hspace{6mm} \text{ for every }
q\in U_{\delta_1,j}^i.
\eeqn Let $G=\{(x,y)\in U_{\delta_{2},j}^i \times \R^{m+1}: y \in \omega(p)\}$.
We note here that $\omega(p)$ does not depend on $x$.
Let the function $F$ be defined by \\
\parbox{14cm}{\beqnn F:\;\;\;\;\; \,G\;\;&\rightarrow& \mathbb{R}^{m+1},\\
(x,y)&\mapsto&f(x)+y, \eeqnn} \hfill
\parbox{8mm}{\beqn \eeqn} \\
where $y\in \omega(p)$. With an argumentation as in Lemma \ref{tubularF},
using (\ref{anglepfixed}) and the fact that $\omega(p)$ is constant,
we conclude that $F(G)$ forms a
tubular neighborhood around $f^i(U_{\delta_2,j}^i)$, and moreover
$B_\sigma(f^i(U_{\delta_{3},j}^i))\subset F(G)$ with $\sigma$ as in (\ref{def3sigma}). \\ \\
We would like to show that
$f^1(U_{\delta_3,j}^1)\subset B_\sigma(f^i(U_{\delta_{3},j}^i))$.
For that let
$p' \in U_{\delta_3,j}^1$. Then there is a unique $x\in B_{\delta_3}$ with
$f^1(p')=A_j^1(x,u_j^1(x))$. Moreover there is a unique $q'\in U_{\delta_3,j}^i$
with $f^i(q')=A_j^i(x,u_j^i(x))$. We estimate
\beqnn
|f^i(q')-f^1(p')|&=&|A_j^i(x,u_j^i(x))-A_j^1(x,u_j^1(x))| \\
&=& |R_j^i(x,u_j^i(x))+T_j^i-R_j^1(x,u_j^1(x))-T_j^1| \\
&\leq& |R_j^i(x,u_j^i(x))-R_j^i(x,u_j^1(x))|+|R_j^i(x,u_j^1(x))-R_j^1(x,u_j^1(x))|
+|T_j^i-T_j^1| \\
&=&|R_j^i((x,u_j^i(x))-(x,u_j^1(x)))|+|(R_j^i-R_j^1)(x,u_j^1(x))|+|T_j^i-T_j^1| \\
&\leq& |u_j^i(x)-u_j^1(x)|+\|R_j^i-R_j^1\||(x,u_j^1(x))|+|T_j^i-T_j^1| \\
&<& \frac{\sigma}{3}+\frac{\sigma}{3}+\frac{\sigma}{3}\\
&=&\sigma, \eeqnn
where in the sixth line we used $|(x,u_j^1(x))|\leq (1+\lambda)r$ and $\mathfrak{d}(\Gamma^1,\Gamma^i)<
[3(1+\lambda)(1+r)]^{-1}\sigma$ which follows from (\ref{dlesssigma}).
Hence
$f^1(U_{\delta_3,j}^1)\subset B_\sigma(f^i(U_{\delta_{3},j}^i))$, i.e.\ $f^1(U_{\delta_3,j}^1)$
lies within the tubular neighborhood defined above.
But this means that there is a $q\in U_{\delta_2,j}^i$ such that $f^1(p)+\omega(p)$
equals $f^i(q)+\omega(p)$. Hence $f^1(p)+\omega(p)$ intersects $f^i(U_{\delta_{2},j}^i)$ in
the point $f^i(q)$. \\ \\
It remains to show that $f^1(p)+\omega(p)$ intersects $f^i(U_{\delta_1,j}^i)$ in not more
than one point. By (\ref{anglepfixed}) we have $\sphericalangle(T(p),\nu_j^i(q))<\frac{\pi}{2}$
for every $q\in U_{\delta_1,j}^i$. By the definition
of $\omega$ this implies $\R^{m+1}=\tau_{f^{i}}(q)\oplus \omega(p)$ for every
$q\in U_{\delta_1,j}^i$. As $f^i$ is an $(r,\lambda)$-immersion, we conclude with
Lemma \ref{intersectionmostonepoint} in the appendix that $f^1(p)+\omega(p)$
intersects $f^i(U_{\delta_1,j}^i)$ in at most one point. \hfill $\square$
\vspace{10mm} \\
\noindent The following lemma will be needed in order to show that the
mappings $\phi^i$ are well-defined: \\[-2mm]
\begin{lemma} \label{welldefined} Let $p\in U_{\delta_{3},j}^1\cap U_{\delta_{3},k}^1$.
Let $S_1$ be the intersection point of
$f^1(p)+\omega(p)$ with $f^i(U_{\delta_{1},j}^i)$, and $S_2$ the intersection point of
$f^1(p)+\omega(p)$ with $f^i(U_{\delta_{1},k}^i)$. Finally let $\sigma_1\in U_{\delta_{1},j}^i$
with $f^i(\sigma_1)=S_1$, and $\sigma_2\in U_{\delta_{1},k}^i$
with $f^i(\sigma_2)=S_2$. Then $\sigma_1=\sigma_2$. \end{lemma}
\noindent \textbf{Proof:} \\
By Lemma \ref{intersectionpoint} we have $S_2\in f^i(U_{\delta_{2},k}^i)$, hence $\sigma_2 \in U_{\delta_{2},k}^i$.
As $p\in U_{\delta_{3},j}^1\cap U_{\delta_{3},k}^1$, we have in particular
$U_{\delta_{2},j}^1\cap U_{\delta_{2},k}^1 \neq \emptyset$ and hence
$U_{\delta_2,k}^1\subset U_{\delta_1,j}^1$ by Lemma \ref{intersect} b). By
Lemma \ref{intersectionpoint} the set $f^1(p)+\omega(p)$
has exactly one point of intersection with $f^i(U_{\delta_{1},j}^i)$. We conclude
$\sigma_1=\sigma_2$. \hfill $\square$ \\
\vspace{6mm} \\
Now we are able to define the mappings $\phi^i:M^1\rightarrow M^i$.
Let $p \in M^1$. Then $p \in U_{\delta_{3},j}^1$ for some $j$. The line
$f^1(p)+\omega(p)$ intersects $f^i(U_{\delta_{1},j}^i)$ in exactly
one point $S_p$. Moreover there is exactly one point
$\sigma_p\in U_{\delta_{1},j}^i$
with $f^i(\sigma_p)=S_p$. We set $\phi^i(p):=\sigma_p$.
The mappings $\phi^i$ are well-defined by Lemma \ref{welldefined}.
Clearly we have $f^i\circ \phi^i(p)=S_p$. \\ \\
We like to show that each $\phi^i$ is a diffeomorphism. For that
we follow in parts the argumentation of \cite{breuning1}: \\
\begin{lemma} \label{surjective} Each of the mappings $\phi^i:M^1\rightarrow M^i$ is surjective. \end{lemma}
\noindent \textbf{Proof:} \\
Let $q\in M^i$. As $Q^i$ is a $\delta_4$-net for $f^i$, there is a $j\in\{1,\ldots,s\}$
with $q\in U_{\delta_4,j}^i$.
By Lemma \ref{tubularF},
for $\varepsilon=\frac{1}{L}\cos \gamma$ the set $F(E_j^\varepsilon)$ forms a tubular
neighborhood around $f^1(U_{\delta_3,j}^1)$, and
moreover $B_\sigma(f^1(U_{\delta_4,j}^1))\subset F_j(E_j^\varepsilon)$ with $\sigma$ as
in (\ref{def3sigma}).
With (\ref{dlesssigma}) and an estimation completely analogous to that in the proof
of Lemma \ref{intersectionpoint}, one shows
$f^i(U_{\delta_4,j}^i)\subset B_\sigma(f^1(U_{\delta_4,j}^1))$.
Hence, for every $q\in U_{\delta_4,j}^i$ there is a $p\in U_{\delta_3,j}^1$ with
$f^i(q)\in f^1(p)+\omega(p)$. By the definition of $\phi^i$ this yields
$\phi^i(p)=q$.
 \hfill
$\square$ \\
\begin{lemma} \label{injective} Each of the mappings $\phi^i:M^1\rightarrow M^i$ is injective. \end{lemma}
\noindent \textbf{Proof:} \\
First we note that for every $j\in \{1,\ldots,s\}$
we have $\phi^i(U_{\delta_5,j}^1)\subset U_{\delta_4,j}^i$. This is shown by the same argumentation
as in Lemma \ref{intersectionpoint}. Moreover, by the proof of Lemma \ref{surjective}, we know that
$f^i(U_{\delta_4,j}^i)\subset F_j(E_j^\varepsilon)$. Using
that $Q^1$ is a $\delta_5$-net for $f^1$,
we conclude $f^i\circ \phi^i(x) \in F_j(\hat{E}_x\cap E_j^\varepsilon)$ for every $x\in U_{\delta_3,j}^1$
(where $\hat{E}_x=\hat{\pi}^{-1}(x)$). As $F_j|E_j^\varepsilon$ is a diffeomorphism, we
conclude that $\phi^i$ is injective on $U_{\delta_3,j}^1$. \\ \\
For showing global injectivity, let $x,y\in M^1$ with $x\neq y$. As $Q^1$ is a
$\delta_5$-net for $f^1$, there are $j,k$ with
$x\in U_{\delta_5,j}^1\subset U_{\delta_4,j}^1$, $y \in U_{\delta_5,k}^1\subset U_{\delta_4,k}^1$. \\ \\
\textbf{Case 1:}\, $U_{\delta_4,j}^1\cap U_{\delta_4,k}^1= \emptyset$ \\ \\
By the considerations at the beginning of this proof,
we have $\phi^i(x)\in U_{\delta_4,j}^i$, $\phi^i(y)\in U_{\delta_4,k}^i$. By (\ref{Z4}) for
$\iota=4$, we have $U_{\delta_4,j}^i\cap U_{\delta_4,k}^i= \emptyset$. This implies
$\phi^i(x)\neq \phi^i(y)$. \\ \\
\textbf{Case 2:}\, $U_{\delta_4,j}^1\cap U_{\delta_4,k}^1\neq \emptyset$ \\ \\
By Lemma \ref{intersect} b) we have $U_{\delta_4,k}^1\subset U_{\delta_3,j}^1$.
By the considerations of above, $\phi^i$ is injective on $U_{\delta_3,j}^1$. Again we conclude $\phi^i(x)\neq \phi^i(y)$.
\hfill $\square$
\vspace{10mm}
\begin{cor} Each mapping $\phi^i:M^1\rightarrow M^i$ is a diffeomorphism. \end{cor}
\noindent \textbf{Proof:} \\
As in Lemma \ref{injective} we have $f^i\circ \phi^i(x) \in F_j(\hat{E}_x\cap E_j^\varepsilon)$ for every $x\in U_{\delta_3,j}^1$.
Using a trivialization of the trivial bundle $\hat{E}_j$, one easily concludes that
$f^i\circ \phi^i:M^1\rightarrow \R^{m+1}$ is an immersion (see also \cite{breuning1}). Moreover, the mapping $\phi^i$ is surjective
by Lemma \ref{surjective}, and injective by Lemma \ref{injective}. We conclude that $\phi^i$ is a diffeomorphism.
\hfill $\square$
\vspace{10mm} \\
Finally we would like to prove that the reparametrizations $f^i\circ \phi^i$ are uniformly Lipschitz bounded.
As above, for $j\in\{1,\ldots,s\}$ we can consider $f^i\circ \phi^i|U_{\delta_3,j}^1$
also as a mapping defined on $B_{\delta_3}$. This mapping shall be denoted by $\hat{f}^i:B_{\delta_3}\rightarrow \R^{m+1}$.
\\
\begin{lemma} \label{Lipschitzrep} Let $j\in\{1,\ldots,s\}$. Let $\hat{f}^i:B_{\delta_{3}}\rightarrow \R^{m+1}$ be the local
representation of $f^i\circ\nolinebreak \phi^i|U_{\delta_{3},j}^1$ as explained above. Then $\hat{f}^i$ is $\Lambda$-Lipschitz
for a finite constant $\Lambda=\Lambda(\lambda)$.
\end{lemma}
\noindent \textbf{Proof:} \\
Let $x,y\in B_{\delta_3}$.
Then there are unique $\mu_1,\mu_2\in \R$ such that \beqnn
\hat{f}^i(x)=(x,u_j^1(x))+\mu_1T(x), \hspace{10mm} \hat{f}^i(y)=(y,u_j^1(y))+\mu_2T(y). \eeqnn
By the construction of the mappings $\phi^i$ we have $|\mu_1|,|\mu_2|<\varepsilon$, where
$\varepsilon=\frac{1}{L}\cos \gamma<r$.
Let $E\in G_{m+1,m}$ be the $m$-space perpendicular to $T(x)$. We define an affine
subspace $\tilde{E}:=(x,u_j^1(x))+E$. Let $\tilde{\pi}:\R^{m+1}\rightarrow \tilde{E}$
denote the orthogonal projection onto $\tilde{E}$. As \beqnn
\tilde{\pi}((x,u_j^1(x))+\mu T(x))=(x,u_j^1(x)) \eeqnn
for any $\mu\in \R$, we may estimate as follows: \beqn \label{Lipschitz1}
|\tilde{\pi}(\hat{f}^i(x))-\tilde{\pi}(\hat{f}^i(y))|
&=& |\tilde{\pi}((x,u_j^1(x))+\mu_1T(x))
-\tilde{\pi}((y,u_j^1(y))+\mu_2T(y))| \nonumber \\
&=& |\tilde{\pi}((x,u_j^1(x))+\mu_2T(x))
-\tilde{\pi}((y,u_j^1(y))+\mu_2T(y))| \nonumber \\
&\leq& |(x,u_j^1(x))-(y,u_j^1(y))+\mu_2(T(x)-T(y))| \\
&\leq& |x-y|+|u_j^1(x)-u_j^1(y)|+r|T(x)-T(y)| \nonumber \\
&\leq& (1+\lambda+rL)|x-y|. \nonumber \eeqn
By Lemma \ref{intersectionmostonepoint} together with Lemma \ref{sestimate} b),
the set $f^i(U_{\delta_{1},j}^i)$ is the graph of a function $\tilde{u}$ on an open subset $U$
of $\tilde{E}$. In the same manner, $f^i(U_{\delta_{2},j}^i)$
is the graph of the same function restricted to a subset $V\subset\!\subset U$.
Again by Lemma \ref{sestimate} b),
on convex subsets of $U$ the function $\tilde{u}$ is $\lambda'$-Lipschitz with
$\lambda'=\tan \gamma$, where $\gamma$ is as in (\ref{def2gamma}).
Let $\varrho>0$ be small enough, such that $B_\varrho(\xi)\subset U$ for any
$\xi\in V$ (where here $B_\varrho(\xi)$ denotes an open ball in $\tilde{E}$). \\ \\
Now assume $|x-y|<\frac{\varrho}{1+\lambda+rL}$.
By Lemma \ref{intersectionpoint} we have $\hat{f}^i(z)\in f^i(U_{\delta_2,j}^i)$ for any $z\in B_{\delta_3}$.
Hence by (\ref{Lipschitz1}) the points
$\tilde{\pi}(\hat{f}^i(x))$ and $\tilde{\pi}(\hat{f}^i(y))$ lie both in the convex subset
$B_\varrho(\tilde{\pi}(\hat{f}^i(x)))$ of $U$. We conclude \\
\parbox{14cm}{\beqnn
|\hat{f}^i(x)-\hat{f}^i(y)|&=& |(\tilde{\pi}(\hat{f}^i(x)),\tilde{u}(\tilde{\pi}(\hat{f}^i(x))))-
(\tilde{\pi}(\hat{f}^i(y)),\tilde{u}(\tilde{\pi}(\hat{f}^i(y))))| \\
&\leq&(1+\tan \gamma)(1+\lambda+rL)|x-y|. \eeqnn} \hfill
\parbox{8mm} {\beqn \label{Lipschitzxysmall} \eeqn} \\
If $x,y\in B_{\delta_3}$ are arbitrary points, let $N\in\N$ with $N>\frac{1+\lambda+rL}{\varrho}|x-y|$.
We define
$x_\iota=x+\iota\,\frac{y-x}{N}\in B_{\delta_3}$ for $\iota=0,\ldots,N$.
Then, using a telescoping sum and (\ref{Lipschitzxysmall}), we have
\beqnn
|\hat{f}^i(x)-\hat{f}^i(y)|&\leq& \sum_{\iota=0}^{N-1}|\hat{f}^i(x_\iota)-\hat{f}^i(x_{\iota+1})| \\
&\leq& (1+\tan \gamma)(1+\lambda+rL)|x-y|. \eeqnn
By the definitions of $L$ and $\gamma$, the quantities $rL$ and $\gamma$ depend only on $\lambda$.
Hence $\hat{f}^i$ is \linebreak
$\Lambda$-Lipschitz with $\Lambda=\Lambda(\lambda)=(1+\tan \gamma)(1+\lambda+rL)$.
\hfill $\square$ \\
\begin{rem} If we choose some of the constants more carefully, we can give
a better bound for $\Lambda$ in the preceding lemma. Choosing the right hand side in (\ref{cauchydist})
extremely small, we can replace $\gamma$ by a number $\tilde{\gamma}$ which is slightly greater than $\arctan \lambda$.
Moreover, we can choose $\varepsilon$ with $|\mu_1|,|\mu_2|<\varepsilon$ so small, that the term
$\varepsilon L$ can almost be neglected.
With these constants, we finally obtain
$\Lambda=(1+ \tan \tilde{\gamma})(1+\lambda+\varepsilon L)<2(1+\lambda)^2$. In particular,
$\Lambda$ does not depend on the dimension $m$ here, although $L$ depends on $m$.
\end{rem}
\noindent Finally, with Lemma \ref{Lipschitzrep},
we may pass to a subsequence such that $f^i\circ \phi^i$ converges uniformly
to a limit function $f:M^1\rightarrow \R^{m+1}$.
As limit manifold we define $M:=M^1$. Thus the limit manifold is a compact
differentiable $m$-manifold.
\end{section}
\begin{section}[The limit function
lies in $\mathfrak{F}^0(r,\lambda)$]{The limit function
lies in \boldmath$\mathfrak{F}^0(r,\lambda)$\unboldmath}
Up to this point we have found a subsequence and diffeomorphisms $\phi^i:M^1\rightarrow M^i$, such
that $f^i\circ \phi^i$ is uniformly Lipschitz bounded and converges uniformly to
an $f:M^1\rightarrow \R^{m+1}$.
In this section we will show that the limit function $f$ lies in $\mathfrak{F}^0(r,\lambda)$. \\ \\
For that we have to show, that
for each point $q\in M^1$ there is an $E=E(q)\in G_{m+1,m}$, such that
$f$ is injective on $U_{r,q}^E$ and the set
 $A_{q,E}^{-1}\circ f(U_{r,q}^E)$ is the graph of a Lipschitz continuous function
 $u:B_r\rightarrow \R$ with Lipschitz constant $\lambda$. \\ \\
So let $q \in M^1$. Let $q^i=\phi^i(q)\in M^i$. As each $f^i$ is an $(r,\lambda)$-immersion,
there are $E^i\in G_{m+1,m}$ such that for each $i$ the set
$(A_{q^i,E^i}^i)^{-1}\circ f^i(U_{r,q^i}^{E^i})$ is the graph
of a differentiable function $u^i:B_r\rightarrow \R$ with $\|Du^i\|_{C^{0}(B_{r})} \leq
\lambda$. \\ \\
Passing to another subsequence, we may assume \beqnn
u^i&\rightarrow& u \hspace{4mm} \text{ uniformly}, \\
E^i&\rightarrow& E \hspace{3.5mm} \text{ for the metric } d
\text{ defined in (\ref{induceddistance})} \eeqnn
as $i\rightarrow\infty$, where $u:B_r\rightarrow \R$ and $E\in G_{m+1,m}$.
In particular, $u$ is Lipschitz continuous with Lipschitz constant $\lambda$. \\ \\
Let $A_{q,E}$ be a Euclidean isometry, which maps the origin to $f(q)$, and the subspace
$\R^m\times\{0\}\subset \R^m\times \R$ onto $f(q)+E$.
Then we have in any case $A_{q,E}(\{(x,u(x)):x \in B_r\})\subset f(M^1)$. \\ \\
To finish the proof, we show that $f$ is injective on $U_{r,q}^E$ and that
$A_{q,E}^{-1}\circ f(U_{r,q}^E)$
is the graph of the function $u$. This is true, if and only if for every $\varrho$ with $0<\varrho<r$
the function $f$ is injective on $U_{\varrho,q}^E$ and the
set $A_{q,E}^{-1}\circ f(U_{\varrho,q}^E)$ is the graph of the function
$u|B_\varrho$. \\ \\
We first show the graph property.
Let a $\varrho$ with $0<\varrho<r$ be given. Let $\varepsilon>0$ with
$\varepsilon<\min\{\varrho,r-\varrho\}$. Moreover, let
$U_\varrho^i\subset M^1$ be the $q$-component of the set
$(\pi\circ A_{q^i,E^i}^{-1}\circ f^i\circ \phi^i)^{-1}(B_\varrho)$. Again,
$U_{\varrho,q}^E\subset M^1$ is the $q$-component of
$(\pi\circ A_{q,E}^{-1}\circ f)^{-1}(B_\varrho)$. By the definition of $U_\varrho^i$, we have
$A_{q^i,E^i}^{-1}\circ f^i\circ \phi^i(U_\varrho^i)=\{(x,u^i(x)):x\in B_\varrho\}$.
As $A_{q^i,E^i}^{-1}\circ f^i \circ \phi^i\rightarrow A_{q,E}^{-1}\circ f$ uniformly, we
conclude with the definitions of $U_\varrho^i$ and $U_{\varrho,q}^E$ that
$U_{\varrho-\varepsilon}^i\subset U_{\varrho,q}^E\subset U_{\varrho+\varepsilon}^i$
for $i$ sufficiently large, in particular
\beqnn
\{(x,u^i(x)):x\in B_{\varrho-\varepsilon}\}\,\subset\, A_{q^i,E^i}^{-1}\circ f^i\circ \phi^i(U_{\varrho,q}^E)
\,\subset\, \{(x,u^i(x)):x\in B_{\varrho+\varepsilon}\}. \eeqnn
Letting $i\rightarrow\infty$, we obtain
\beqnn
\{(x,u(x)):x\in B_{\varrho-\varepsilon}\}\,\subset\, A_{q,E}^{-1}\circ f(U_{\varrho,q}^E)
\,\subset\, \{(x,u(x)):x\in B_{\varrho+\varepsilon}\}. \eeqnn
As this is true for every $\varepsilon>0$ with $\varepsilon<\min\{\varrho,r-\varrho\}$,
we conclude with the definition of $U_{\varrho,q}^E$ that
$A_{q,E}^{-1}\circ f(U_{\varrho,q}^E)=\{(x,u(x)):x\in B_\varrho\}$.
This is the desired graph property. \\ \\
Similarly, one shows that $f$ is injective on $U_{\varrho,q}^E$. We have
$f(x)=\lim_{i\rightarrow\infty}f^i\circ \phi^i(x)$ for all $x\in U_{\varrho,q}^E$,
and moreover $U_{\varrho,q}^E\subset U_{\varrho+\varepsilon}^i$ for $i$ sufficiently large.
The functions $f^i\circ \phi^i$ are injective on $U_{\varrho+\varepsilon}^i$ and it holds
$A_{q^i,E^i}^{-1}\circ f^i\circ \phi^i(U_{\varrho+\varepsilon}^i)=
\{(x,u^i(x)):x\in B_{\varrho+\varepsilon}\}$. Using $A_{q^i,E^i}\rightarrow A_{q,E}$, one
easily concludes that $A_{q,E}^{-1}\circ f$ and hence also $f$ is injective on $U_{\varrho,q}^E$. \\ \\
This shows that the limit function $f$ lies in $\mathfrak{F}^0(r,\lambda)$.
\end{section}
\begin{section}{Compactness in higher codimension} \label{compactnesshigher}
In the final section we want to prove Theorem \ref{compactness3},
that is compactness of $(r,\lambda)$-immersions in higher codimension
with $\lambda\leq\frac{1}{4}$.
Our main task here is to give an analogous construction of the averaged normal projection for arbitrary
codimension. For that we shall use a Riemannian center of mass,
which was introduced above. \\ \\
So let $f^i$ be a sequence as in Theorem \ref{compactness3} with $\lambda\leq \frac{1}{4}$.
For all objects of the preceding sections that are defined also in arbitrary codimension,
we shall use precisely the same notation. We note that
Lemmas \ref{intersect} and \ref{cover} are true also in higher codimension.
For $q\in M^1$ we set \beqnn
\lambda_j^q:=g\left(\hspace{0.2mm}\frac{|f^1(q)-f^1(q_j)|}{\delta_2}\right). \eeqnn
As in the proof of Lemma \ref{sestimate} a) we conclude that there is a $k\in Z(q)$ with
$\lambda_k^q=1$.
For each $j\in \{1,\ldots,s\}$ let $N_j\in G_{n,k}$ be the
$k$-space perpendicular to $E_j$.
We define for each $q\in M^1$ a probability measure $\mu_q$ on $G_{n,k}$ by \beqnn
\mu_q=\Biggl(\sum_{j\in Z(q)}\lambda_j^q\Biggr)^{\!-1}\!\!\sum_{j\in Z(q)}\lambda_j^q\delta_{N_j}\,, \eeqnn
where $\delta_N$ denotes the Dirac measure on $G_{n,k}$ supported at $N\in G_{n,k}$. \\ \\
Moreover, let \beqnn
\nu:M^1&\rightarrow& G_{n,k}, \\
q&\mapsto& (T_qM^1)^\bot \eeqnn
be the normal map of $f^1$ as defined in (\ref{normalnotion}),
and $\tau:M^1\rightarrow G_{n,m}$ the corresponding tangent map as in (\ref{notiontangentspace}).
Now consider \beqnn
P: \overline{B}_{\frac{\pi}{6}}(\nu(q))&\rightarrow& \R, \\
P(p)&=&\int_{G_{n,k}} d(p,x)^2\,d\mu_q(x), \eeqnn
where $\overline{B}_{\frac{\pi}{6}}(\nu(q))\subset G_{n,k}$ is the closed ball of radius
$\frac{\pi}{6}$ around $\nu(q)$. Here the radius is measured with respect
to the canonical distance $d$ on $G_{n,k}$ as defined
in (\ref{induceddistance}). \\
\begin{lemma} \label{sptpi12}
For every $q\in M^1$
it holds\, $\spt \mu_q\subset B_{\frac{\pi}{12}}(\nu(q))$. \end{lemma}
\textbf{Proof:} \\
By the definition of $\mu_q$ it is sufficient to show that $N_j$ lies in
$B_{\frac{\pi}{12}}(\nu(q))$ for every $j\in Z(q)$.
So let $j\in Z(q)$. By the definition of $Z(q)$ we have $q\in U_{\delta_2,j}^1$.
We deduce that $N_j$ is the graph of a linear function $h$ over $\nu(q)$ with
$\|Dh\|=(\sum_{i=1}^k |\partial_ih|^2)^{\frac{1}{2}}\leq \lambda \leq \frac{1}{4}$.
Let $\theta_1,\ldots,\theta_k$ be the principal angles
between $N_j$ and $\nu(q)$.
After a suitable rotation we may assume that $\tan \theta_i=|\partial_i h|$ for every
$i\in \{1,\ldots,k\}$. Using
$\theta \leq \tan \theta$ for $\theta\in [0,\frac{\pi}{2})$, we estimate
$d(N_j,\nu(q))=(\sum_{i=1}^k\theta_i^2)^{\frac{1}{2}}
\leq(\sum_{i=1}^k(\tan \theta_i)^2)^{\frac{1}{2}}=(\sum_{i=1}^k |\partial_ih|^2)^{\frac{1}{2}}
\leq \lambda\leq \frac{1}{4}<\frac{\pi}{12}$. Hence
$N_j$ lies in $B_{\frac{\pi}{12}}(\nu(q))$. \hfill $\square$ \\ \\ \\
In particular we have $\spt \mu_q\subset B_{\frac{\pi}{6}}(\nu(q))$.
Hence we conclude with Lemma \ref{sectcurgrass}
and Theorem \ref{centerofmass}, that there is exactly one
center of mass $N(q)\in B_{\frac{\pi}{6}}(\nu(q))\subset G_{n,k}$ for $\mu_q$.
In this way we may define a mapping \beqnn
N:M^1&\rightarrow& G_{n,k}, \\
q&\mapsto& N(q). \eeqnn
An important property of the averaged normal $N$ constructed in this way is its differentiability.
It is needed in order to obtain diffeomorphisms $\phi^i:M^1\rightarrow M^i$.
We will show
that $N$ is in $C^\K$ if the function $f^1$ is in $C^\K$
(here we denote by $\K$ the degree of differentiability, and by $k$ the codimension).
First, for functions defined on manifolds, we need
the following variation of the implicit function theorem: \\
\begin{lemma} \label{implicitformanifolds}
Let $M$ be a smooth $m$-manifold, $(N,g)$ a smooth Riemannian $n$-manifold and
$f:M\times \nolinebreak N\rightarrow \R$ a mapping. For every fixed $x\in M$, assume that \vspace{-1mm}
\beqnn h_x:N\rightarrow \R,\hspace{4mm} h_x=f(x,\cdot\,) \eeqnn \vspace{-6mm} \\
is in $C^2(N)$ and is strictly convex. Let \,$\K\geq 1$ be an integer. Denoting by\,
$\grad h_x$ the gradient of the fixed function $h_x$ defined above,
assume that \vspace{-1mm}
\beqnn H:M\times N\rightarrow TN,\hspace{4mm} (x,y)\mapsto \grad h_x(y) \eeqnn \vspace{-6mm} \\
is in $C^\K(M\times N,TN)$.
Let $(x_0,y_0)\in M\times N$ be a point with $H(x_0,y_0)=0 \in T_{y_0}N$. \\ \\[-2mm]
Then there are open neighborhoods $U\subset M$ of \,$x_0$ and $V\subset N$ of \,$y_0$,
and moreover a function $F\in C^\K(U,V)$, such that
$\{(x,y)\in U\times V: H(x,y)=0\in T_yN\}\;=\;\{(x,F(x)):x \in U\}$.
\end{lemma}
\textbf{Proof:} \\
Let $\varphi_1:U_1\rightarrow \varphi(U_1)$ be a coordinate chart of $M$
with $x_0\in U_1$, and let $\varphi_2:V_1\rightarrow \varphi_2(V_1)$ a coordinate
chart of $N$ with $y_0\in V_1$. For fixed $x\in M$, in the local coordinates
$\varphi_2$ we have  \beqn \label{gradlocalrepres}
\grad h_x =\sum_{i,j=1}^ng^{ij}\partial_jh_x\hspace{1pt}\partial_i, \eeqn
and, with the corresponding Christoffel symbols $\Gamma_{ij}^k=\frac{1}{2}\sum_{l=1}^ng^{kl}(\partial_ig_{jl}+\partial_jg_{il}
-\partial_lg_{ij})$, the components of the Hessian
$D_{ij}^2h_x=\partial_i\partial_j h_x-\sum_{k=1}^n\,\Gamma_{ij}^k\, \partial_k h_x$.
If we assume $\varphi_2$ to be Riemannian normal coordinates centered in $y_0$, we obtain
\beqn \label{Hessianmatrix}
D_{ij}^2h_x(y_0) = \partial_i\partial_jh_x(y_0). \eeqn
Let us now consider the local representations of $h_x$ and $f$ in the coordinates
$\varphi_2$ and $\varphi_1\times \varphi_2$ respectively. We denote these representations
simply by $h_x$ and $f$ again. Moreover, we identify $x_0$ and $y_0$ with $\varphi_1(x_0)$ and $\varphi_2(y_0)$
respectively. The condition on $h_{x}$ to be strictly convex means
that the Hessian $D^2h_x$ is positive definite in every point.
Hence, by (\ref{Hessianmatrix}), the Hessian matrix $D^2h_x(y_0)$ of the local representation
is positive definite, in particular \beqn \label{Hessianinvertible}
D^2h_{x_0}(y_0) \hspace{3mm} \text{ is invertible.} \eeqn
The Jacobian $Df$ may be considered as a mapping $Df\!:\Omega\rightarrow \R^{m+n}$,
where $\Omega\!=\!\varphi_1(U_1)\times\nolinebreak \varphi_2(V_1)\linebreak\subset \R^m\times \R^n$.
We write $Df=(D_xf,D_yf)\in \R^m\times \R^n$ and consider the mapping
$D_yf:\Omega\rightarrow \R^n$. Similarly, for the Jacobian of $D_yf$, we write
$D(D_yf)=(D_x(D_yf),D_y(D_yf))\in \R^{n\times m}\times \R^{n\times n}$.
As $D_yf(x_0,y_0)=Dh_{x_0}(y_0)$ and as $H(x_0,y_0)=0$, we conclude \beqn \label{implicitcond1}
D_yf(x_0,y_0)=0. \eeqn
Similarly, as $D_y(D_yf)(x_0,y_0)=D^2h_{x_0}(y_0)$, we know by (\ref{Hessianinvertible}) that
\beqn \label{implicitcond2}
D_y(D_yf)(x_0,y_0) \hspace{3mm} \text{ is invertible.} \eeqn
The assumption on $H$ to be in $C^\K$ implies with (\ref{gradlocalrepres})
that also $D_yf:\Omega\rightarrow \R^n$ is
in $C^\K$. Hence we may use (\ref{implicitcond1}), (\ref{implicitcond2}) and apply the usual implicit function
theorem to the function $D_yf$. From this we deduce the statement. \hfill $\square$ \\ \\
Using the preceding lemma, we are able to deduce that the mapping $N$ is differentiable: \\[-2mm]
\begin{lemma} Let $N:M^1\rightarrow G_{n,k}$ be the
averaged normal corresponding to $f^1:M^1\rightarrow \R^n$,
as constructed above. Assume that $f^1\in C^\K(M^1,\R^n)$ for a $\K\geq 1$. Then
$N\in C^\K(M^1,G_{n,k})$. \end{lemma}
\textbf{Proof:} \\
Let $q_{_0}\in M^1$ be a point. We show that $N$ is $C^\K$ in a neighborhood of $q_{_0}$.
Let $W\subset M^1$ be an open neighborhood of $q_{_0}$ with $\nu(W)\subset B_{\frac{\pi}{12}}(\nu(q_{_0}))$.
With Lemma \ref{sptpi12} we have $\spt \mu_q\subset B_{\frac{\pi}{6}}(\nu(q_{_0}))$ for every $q\in W$;
this will be implicitly used in the following argumentation.
Let \beqnn
G:  W \times B_{\frac{\pi}{6}}(\nu(q_{_0})) &\rightarrow& \R, \\
(q,p)&\mapsto& \int_{G_{n,k}} d(p,x)^2\,d\mu_q(x). \eeqnn Moreover, for fixed $q \in W$ let
$h_q:B_{\frac{\pi}{6}}(\nu(q_{_0}))\rightarrow \R$, $h_q:=G(q,\cdot\,)$.
By this definition, $h_q$ is smooth on $B_{\frac{\pi}{6}}(\nu(q_{_0}))$ and by Theorem \ref{centerofmass} strictly
convex. We denote by $\grad h_q$ the gradient of the fixed function $h_q$, and define \beqnn
H:W\times B_{\frac{\pi}{6}}(\nu(q_{_0}))\rightarrow TB_{\frac{\pi}{6}}(\nu(q_{_0})),
\hspace{4mm} (q,p)\mapsto \grad h_q(p). \eeqnn
With (\ref{gradexpformula}) and the definition of $\mu_q$, we calculate \beqn \label{Hrepresentation}
H(q,p)=-2\Biggl(\sum_{j\in Z(q)}\lambda_j^q\Biggr)^{\!-1}\!\!\sum_{j\in Z(q)}\lambda_j^q
\exp_p^{-1}(N_j). \eeqn
As $\lambda_j^q=g\left(\frac{|f^1(q)-f^1(q_j)|}{\delta_2}\right)$ and by the definition
of $g$, the mapping $q\mapsto \lambda_j^q$ is in $C^\K$ if $f$ is in $C^\K$.
Moreover, as for every $j\in Z(q)$ the mapping $p\mapsto \exp_p^{-1}(N_j)$ is smooth, we conclude
that $H$ is in $C^\K$. Note that $g$ is smooth with $g(1)=0$, hence $H$ is $C^\K$ even if
the sums in (\ref{Hrepresentation}) depend on $Z(q)$. \\ \\
As $N(q)\in B_{\frac{\pi}{6}}(\nu(q_{_0}))$ is the center of mass for $\mu_q$, we have
$H(q,N(q))=0$ for every $q\in W$, in particular
$H(q_{_0},N(q_{_0}))=0$. \\ \\
Now we are in a position to apply Lemma \ref{implicitformanifolds}.
We conclude that there are open neighborhoods $U\subset W$ of
$q_{_0}$, $V\subset B_{\frac{\pi}{6}}(\nu(q_{_0}))$ of $N(q_{_0})$, and a function
$F\in C^\K(U,V)$ with $\{(x,y)\in U\times V:H(x,y)=0\}=\{(x,F(x)):x\in U\}$.
With Theorem \ref{centerofmass} we deduce, that $N$ coincides with $F$ on $U$. Hence $N$ is in
$C^\K$ on $U$. \hfill $\square$
\vspace{3mm}
\begin{rem} \label{reasonaveragednormal}In particular, the preceding lemma shows that the averaged normal $N$ can be used
for the projection in the case of immersions with $L^p$-bounded second fundamental form,
which was the case considered in \cite{breuning1}. For an $(r,\lambda)$-immersion $f\in C^\K$,
the normal $\nu_f$ is in $C^{\K-1}$, while the averaged normal $N$ is in $C^\K$. In particular,
the averaged normal of a $C^1$-immersion is differentiable and forms locally a tubular neighborhood
around the immersion. Thus it is possible to construct diffeomorphisms $\phi^i:M^1\rightarrow M^i$ using the
averaged normal. However, if one likes to show convergence as in \cite{langer} and in \cite{breuning1}, we require
$N$ even to be in $C^{2}$. For that purpose,
an additional smoothing of $f$ is unavoidable; this was also performed by Langer (see the first paragraph
on p.\ 229 in \cite{langer}, where a $C^1$-perturbation is made in order to smooth the immersion).
On the other hand, a pure smoothing argument would not suffice to prove Theorems
\ref{compactness2} and \ref{compactness3}. As in general the limit is not even differentiable, one has
to project from $f^{i_{0}}(M^{i_{0}})$ for a fixed and sufficiently large  $i_0$. The averaged normal is needed
then in order to estimate the size of the tubular neighborhood.
\end{rem}
\vspace{3mm}
As in the case of codimension $1$, we may consider the restriction of $N$ to $U_{\delta_3,j}^1$
as a mapping defined on $B_{\delta_3}$.
As an analogue of Lemma \ref{TisLipschitz} we show the following statement: \\[-2mm]
\begin{lemma}
If we consider $G_{n,k}$ as a metric space with the geodesic distance $d$,
the mapping $N:B_{\delta_3}\rightarrow G_{n,k}$ is $L$-Lipschitz
with $L=4^{12m+6}r^{-1}$. \end{lemma}
\textbf{Proof:} \\
Let $x,y\in B_{\delta_{3}}$. Then there are unique $p,q\in U_{\delta_{3},j}^1$ with
$\pi\circ A_j^{-1}\circ f^1(p)=x$, \;$\pi\circ A_j^{-1}\circ f^1(q)=y$.
With the argumentation at the beginning of the proof of
Lemma \ref{TisLipschitz}, one shows $\lambda_j^p=0$ for $j\in Z(q)\setminus Z(p)$
and $\lambda_j^q=0$ for $j\in Z(p)\setminus Z(q)$. Again as in Lemma \ref{TisLipschitz},
we estimate \\ \parbox{14cm}{\beqnn
|\lambda_j^p-\lambda_j^q|&\leq& 36(1+\lambda)^3r^{-1}|x-y| \\
&\leq&72r^{-1}|x-y|, \eeqnn} \hfill
\parbox{8mm}{\beqn \label{lipschi1} \eeqn} \\
where we used $\lambda\leq\frac{1}{4}$.
Now we note that $\sum_{k\in Z(p)\cup Z(q)}\lambda_k^p \geq 1$ and $\sum_{k\in Z(p)\cup Z(q)}\lambda_k^q\geq 1$.
Moreover $|Z(p)\cup Z(q)|\leq 2[3(1+\lambda)]^{6m}\leq 2\cdot4^{6m}$
for $\lambda\leq \frac{1}{4}$ by (\ref{cardZ}). Using all this, we obtain \\
\parbox{14.6cm} {\beqnn
\left\lvert\Biggl(\hspace{1pt}\sum_{k\in Z(p)}\lambda_k^p\Biggr)^{\!-1}
- \Biggl(\hspace{1pt}\sum_{k\in Z(q)}\lambda_k^q\Biggr)^{\!-1}\!\right\rvert
&\leq& \sum_{k\in Z(p)\cup Z(q)}|\lambda_k^p-\lambda_k^q| \\
&\leq& 9\cdot 4^{6m+2}r^{-1}|x-y|.
\eeqnn} \parbox{8mm}{\beqn \label{lipschi2} \eeqn} \\
Using (\ref{lipschi1}) and (\ref{lipschi2}), one easily concludes \beqn \label{lipschi3}
\left\lvert\Biggl(\hspace{1pt}\sum_{k\in Z(p)}\lambda_k^p\Biggr)^{\!-1}\! \lambda_j^p
- \Biggl(\hspace{1pt}\sum_{k\in Z(q)}\lambda_k^q\Biggr)^{\!-1}\!\lambda_j^q\right\rvert
\leq 10\cdot 4^{6m+2}r^{-1}|x-y|.
\eeqn
Now assume $k\in Z(p)$. Then $U_{\delta_3,j}^1\cap U_{\delta_2,k}^1\neq \emptyset$, hence
$U_{\delta_3,j}^1\subset U_{\delta_1,k}^1$ by Lemma \ref{intersect} b). This implies
$q\in U_{\delta_1,k}^1$. With a calculation as in Lemma \ref{sptpi12} we deduce
$N_k\in B_{\frac{\pi}{12}}(\nu(q))$. We conclude that both $\spt \mu_p$ and $\spt \mu_q$
are a subset of $B_{\frac{\pi}{12}}(\nu(q))$. This enables us to apply Lemma \ref{distcenterofmass}
with $\mu_1=\mu_p$ and $\mu_2=\mu_q$. \\ \\
With Lemma \ref{distcenterofmass}, the definitions of $\mu_p$ and $\mu_q$, and (\ref{lipschi3})
we estimate \beqnn
d(N(x),N(y))&\leq& C\int_{G_{n,k}} d(N(q),z)\,d|\mu_p-\mu_q|(z) \\
&=& C \sum_{j\in Z(p)\cup Z(q)}d(N(q),N_j)\left\lvert\Biggl(\hspace{1pt}\sum_{k\in Z(p)}\lambda_k^p\Biggr)^{\!-1}\! \lambda_j^p
- \Biggl(\hspace{1pt}\sum_{k\in Z(q)}\lambda_k^q\Biggr)^{\!-1}\!\lambda_j^q\right\rvert \\
&\leq& C\cdot 10\cdot 4^{6m+2}\!\!\sum_{j\in Z(p)\cup Z(q)}\!d(N(q),N_j)\,r^{-1}|x-y| \\
&\leq& 4^{12m+6}r^{-1}|x-y|, \eeqnn
where in the last line we used $d(N(q),N_j)<\frac{\pi}{6}$, $|Z(p)\cup Z(q)|\leq 2\cdot4^{6m}$
and $C=\linebreak1+(\kappa^{1/2}\varrho)^{-1}\tan(2 \kappa^{\,1/2}\varrho)
<16$ for $\kappa=2$ and $\varrho=\frac{\pi}{6}$ by (\ref{Ckappar}). \hfill $\square$
\\ \\ \\
Now the rest of the proof is analogous to the case of codimension $1$. First we note that with the preceding
lemma one easily derives an estimate for the size of the tubular neighborhood around $f^1$ formed by
$N$. This is done by using elementary geometry in much the same way as in the appendix (where the case
of codimension $1$ is considered); as we assumed $\lambda$ to be small and hence $N$ nearly to be  perpendicular
to $f^1$, it is even easier here as we can estimate rather roughly (and do not need an analogue of
Lemma \ref{sestimate} b) for that).
Moreover, we can show the existence and uniqueness of intersection points of $f^1(p)+N(p)$ with an
appropriate restriction of $f^i(M^i)$ by the fixed point argument of \cite{breuning1}. To show surjectivity
of $\phi^i$ one uses the estimate for the size of the tubular neighborhood
and shows that $f^i(M^i)$ lies within this neighborhood. The rest of the proof is the same as in the
case of codimension $1$. \\ \\
The question arises, whether compactness in higher codimension, that is Theorem \ref{compactness3},
can also be shown for an arbitrary Lipschitz constant $\lambda$ (as in the case of codimension $1$). Surely,
the bound $\lambda\leq \frac{1}{4}$ is not optimal. One could try to find the largest possible
bound for $\lambda$, and --- in the case that it is finite --- to give a counterexample for immersions exceeding this bound.
We would like to suggest two possibilities for extending the construction in this section
to immersions with Lipschitz constant larger than the ones considered here:
First, as proposed in the remark on p.\ 511 in \cite{karcher}, one could use another definition for the
center of mass, which allows one to define centers in larger balls.
The second is to find a center of mass not in a convex ball, but in a larger convex
subset of $G_{n,k}$. Such kind of subsets of Grassmannians have been detected by J.\ Jost and Y.L.\ Xin in \cite{jxin}.
\end{section} \\
\begin{appendix}
\section{Size of tubular neighborhoods}
In this appendix we like to prove Lemma \ref{tubularsigma}, that is we
estimate the size of a tubular neighborhood around a graph
depending on different quantities such as angles and Lipschitz constants. We shall use the notations
introduced in the paragraph preceding Lemma \ref{tubularsigma}.
For a general treatise on the existence of tubular neighborhoods see \cite{broecker} and \cite{hirsch}.
Moreover, in Lemma \ref{intersectionmostonepoint} we will show a result needed for proving
that the projection in Section 4 has at most one point of intersection with an appropriate
subset of $f^i(M^i)$. \\ \\
\noindent \textbf{Proof of Lemma \ref{tubularsigma}:}
\begin{itemize}
\item[a)] We like to start with the following \emph{initial consideration:} \\ \\
Let $q \in B_\varrho$. Let $f(x)=(x,u(x))$ and $\tau_f(q) \in G_{m+1,m}$ be
the tangent space at $q$ as in (\ref{notiontangentspace}). In particular
$\tau_f(q)$ is an $m$-space in $\R^{m+1}$ perpendicular
to $\nu(q)$. Furthermore let $K\subset \tau_f(q)$ be a $1$-dimensional subspace
of $\tau_f(q)$. Let $p\in B_\varrho$
and let $\alpha\leq\frac{\pi}{2}$ be the smaller angle enclosed by the lines
$\omega(p)$ and $K$.
From (\ref{angelTnu}) we deduce \beqn \label{anglelineest}
\alpha \;\geq\; \frac{\pi}{2}-\gamma\;>\;0. \eeqn
Now let us come to the \emph{main part} of the proof: \\
Let $x,y \in B_\varrho$ with $x\neq y$.
Without loss of generality we may assume $x-y\in \R^1\times\{0\}$\linebreak$\subset \R^m$.
By the mean
value theorem there is a $z\in\{(1-t)x+ty:\:t\in(0,1)\}\subset B_\varrho$ with \beqnn
\partial_1u(z)=\frac{u(x)-u(y)}{x_1-y_1}, \eeqnn
where $x_1, y_1$ are the first coordinate of the vectors $x,y$ respectively.
Let $\{e_1,\ldots,e_m\}$ be the standard basis of $\R^m$. We set \beqnn
K:=\text{span } \{(e_1,\partial_1u(z))\}\subset \tau_f(z). \eeqnn
Let $\alpha\leq\frac{\pi}{2}$ be the smaller angle enclosed by the
lines $\omega(y)$ and $K$.
By (\ref{anglelineest}) we have $\alpha\geq\frac{\pi}{2}-\gamma$.
In particular the smaller angle between $\omega(y)$ and the line through $(x,u(x))$ and $(y,u(y))$
is greater than or equal to $\frac{\pi}{2}-\gamma$ (see Figure \ref{calcdistfig}). \vspace{-0.4cm} \\
\setlength{\unitlength}{2cm}
\begin{picture}(12,8.5)
\put(3,5){\line(1,0){3}}
\put(3,5){\line(-1,0){3}}
\linethickness{0.29mm}
\qbezier(0,5.8)(1,5)(2,5)
\qbezier(2,5)(3,5)(6,5.6)
\thinlines
\put(-0.03,4.96){\small $($}
\put(5.96,4.96){\small $)$}
\put(6.1,4.9){$B_\varrho\cap(\R^1\times\{0\})$}
\put(0.24,4.97){\tiny$|$}
\put(3.21,4.97){\tiny$|$}
\put(0.24,4.8){$x$}
\put(3.21,4.8){$y$}
\put(0.23,5.575){\tiny$\bullet$}
\put(3.2,5.09){\tiny$\bullet$}

\put(0,5.65){\line(6,-1){4.8}}

\put(3.26,5.1){\line(-5,6){2.3}}
\put(3.26,5.1){\line(5,-6){0.7}}
\qbezier[60](0.24,5.59)(1.0023,6.22525)(1.7647,6.8605)
\put(1.75,6.85){\tiny $\bullet$}
\put(1.9,6.9){\small $\pi^\bot((x,u(x)))$}

\put(1.5,7.9){\vector(-1,-4){0.12}}
\put(1.53,7.9){\small $(y,u(y))+\omega(y)$}

\put(4.91,5.73){\small $\{(x,u(x)):x\in B_\varrho\cap(\R^1\times\{0\})\}$}

\put(4.8,4.5){\vector(-3,4){0.25}}
\put(4.8,4.35){\small line through $(x,u(x))$ and $(y,u(y))$}

\put(3.27,5.78){\vector(-1,-2){0.199}}
\put(2.76,5.85){\small an angle $\geq \frac{\pi}{2}-\gamma$}

\qbezier(2.75,5.2)(2.75,5.37)(2.95,5.46)
\put(2.9,5.23){$\alpha$}

\end{picture}
\vspace{-9.5cm} \\
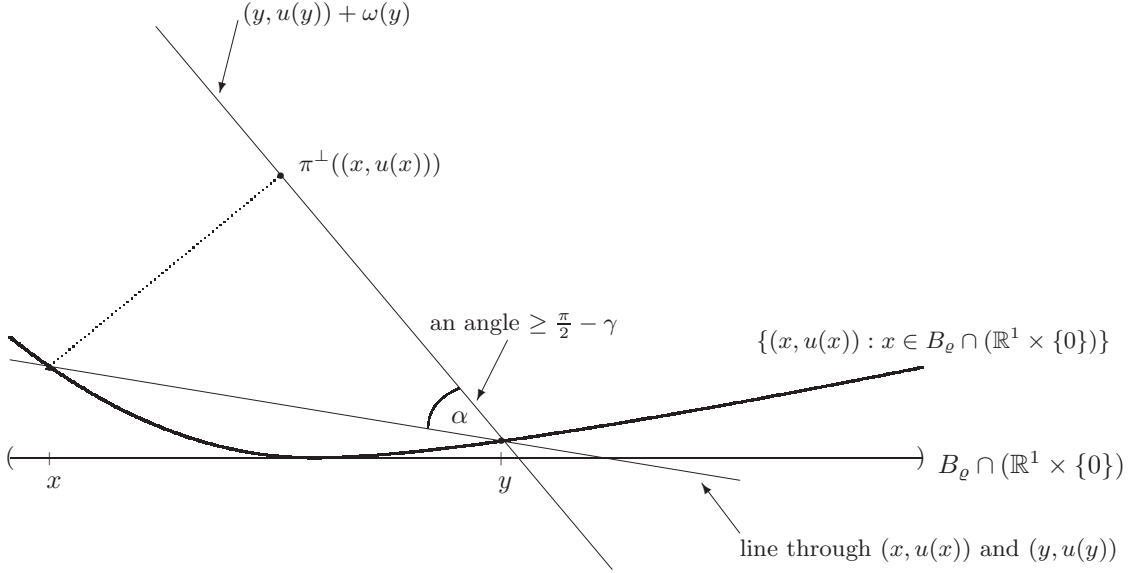
\begin{fig} \label{calcdistfig} Calculation of the distance between $(x,u(x))$ and
$\pi^\bot((x,u(x)))$. Note that, unlike the rest of the figure, the line $(y,u(y))+\omega(y)$ does not necessarily
lie in the plane $(\R^1\times\{0\})\times \R^1\subset \R^{m+1}$.
\end{fig}
\vspace{8mm}
Let $\pi^\bot((x,u(x)))$ denote the orthogonal projection of $(x,u(x))$ onto
$F(\{y\}\times \omega(y))=(y,u(y))+\omega(y)$. Then \beqn
|(x,u(x))-\pi^\bot((x,u(x)))|&\geq&|(x,u(x))-(y,u(y))|\sin\left(\frac{\pi}{2}-\gamma\right) \nonumber \\
&\geq& |x-y|\sin\left(\frac{\pi}{2}-\gamma\right) \label{appendix1} \\
&=& |x-y|\cos \gamma. \nonumber \eeqn \vspace{5mm} \\
Now we distinguish two cases: \\ \begin{itemize}
\item[\textbf{Case 1:}] \beqn [(x,u(x))+\omega(x)]\cap [(y,u(y))+\omega(y)]=\emptyset \label{appendixcase1}\eeqn \\
In this case we do not need any further estimations. \\[2mm]
\item[\textbf{Case 2:}] \beqn [(x,u(x))+\omega(x)]\cap [(y,u(y))+\omega(y)]\neq\emptyset \label{appendixcase2} \eeqn
\vspace{0mm}\\
We now have to consider the following two subcases 2.\textit{i} and 2.\textit{ii:} \\[-3mm]
\begin{itemize}
\item[\textbf{2.\textit{i:}}] The case $|x-y|\leq\frac{1}{L}$. \\ \\
Let $\theta=\sphericalangle(T(x),T(y))$.
By the assumption (\ref{appendixcase2})
we have $\theta>0$.
Using $|T(x)|=|T(y)|=1$, we estimate \beqn
\theta&=&2\arcsin\left(\frac{|T(x)-T(y)|}{2}\right) \nonumber \\
&\leq& 2\arcsin\left(\frac{L}{2}|x-y|\right)\label{appendix2} \\
&<& \frac{\pi}{2}. \nonumber \eeqn
Now let $\xi\in \R^{m+1}$ denote the intersection point of
$(x,u(x))+\omega(x)$ with $(y,u(y))+\omega(y)$. \\\vspace{15mm} \\
\setlength{\unitlength}{1cm}
\begin{picture}(12,13)
\put(6,8){\line(1,0){3}}
\put(6,8){\line(-1,0){3}}
\linethickness{0.29mm}
\qbezier(3,8.8)(4,8)(5,8)
\qbezier(5,8)(6,8)(9,8.6)
\thinlines
\put(2.97,7.92){\footnotesize $($}
\put(8.93,7.92){\footnotesize $)$}

\put(3.24,7.97){\tiny$|$}
\put(6.21,7.97){\tiny$|$}
\put(3.2,7.7){\small $x$}
\put(6.16,7.7){\small $y$}
\put(3.17,8.56){\tiny$\bullet$}
\put(6.2,8.09){\tiny$\bullet$}

\put(3,8.65){\line(6,-1){4.8}}

\put(6.26,8.1){\line(-5,6){5.1}}
\put(6.26,8.1){\line(5,-6){0.7}}
\qbezier[30](3.24,8.59)(4.0023,9.22525)(4.7647,9.8605)
\put(4.75,9.85){\tiny $\bullet$}
\put(4.9,9.9){\small $\pi^\bot((x,u(x)))$}

\put(4.5,10.9){\vector(-1,-4){0.12}}
\put(4.53,10.9){\small $(y,u(y))+\omega(y)$}

\put(3.24,8.575){\line(-1,3){1.9}}

\put(1.44,13.76){\tiny$\bullet$}
\put(1.64,13.8){\small $\xi$}

\qbezier(1.85,12.7)(2.07,12.7)(2.2,12.93)
\put(1.86,12.85){\small $\theta$}

\end{picture}
\vspace{-7cm}
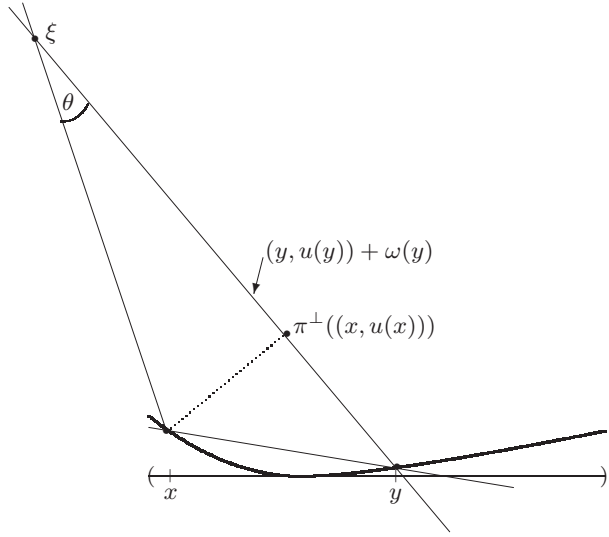
\begin{fig} \label{calcdistfi} Calculation of the distance between $(x,u(x))$ and
$\xi$. Again, \linebreak $(y,u(y))+\omega(y)$ does not necessarily lie in
$(\R^1\times \{0\})\times \R^1$.
\end{fig}
\vspace{1cm}
Then, using (\ref{appendix1}) and (\ref{appendix2}), \beqnn
|(x,u(x))-\xi|&=& \frac{|(x,u(x))-\pi^\bot((x,u(x)))|}{\sin \theta} \\
&\geq& \frac{|x-y|\cos \gamma}{\sin\left(2\arcsin\left(\frac{L}{2}|x-y|\right)\right)} \\
&=& \frac{1}{L} \frac{\cos \gamma}{\sqrt{1-\frac{L^2}{4}|x-y|^2}} \\
&>&\frac{1}{L}\cos \gamma.
\eeqnn
\item[\textbf{2.\textit{ii:}}] The case $|x-y|>\frac{1}{L}$. \\ \\
Let $\xi$ be as in Case 2.\textit{i}. Then (\ref{appendix1}) directly implies \beqnn
|(x,u(x))-\xi|>\frac{1}{L}\cos \gamma. \eeqnn
\end{itemize}
\end{itemize}
Let $\varepsilon=\frac{1}{L}\cos \gamma$. Summarizing Case 1, Case 2.\textit{i} and 2.\textit{ii}, we conclude that
$F$ is injective on $E^\varepsilon$.
Applying well-known results from elementary differential topology, we deduce
that $F|E^\varepsilon$ is a diffeomorphism onto an open neighborhood of
$\{(x,u(x))\in \R^m\times\R:x\in \nolinebreak B_\varrho\}$. This proves part a) of Lemma \ref{tubularsigma}.
\vspace{1.5cm}
\item[b)]
Let $\partial (F(E^\varepsilon))$ denote the boundary of $F(E^\varepsilon)$ in $\R^{m+1}$,
where $\varepsilon=\frac{1}{L}\cos \gamma$ as in part a). Let $x\in\overline{B}_{\frac{\varrho}{2}}$.
We have to show \beqnn
\text{dist}((x,u(x)),\partial (F(E^\varepsilon)))\;\geq\; \sigma \eeqnn
with $\sigma=\min\{\frac{\varrho}{2}\cos \gamma,\frac{\cos^2 \gamma}{2L(1+\lambda)}\}$
as in Lemma \ref{tubularsigma} b). \\ \\ \\
So let $\zeta\in \partial (F(E^\varepsilon))\subset \R^{m+1}$. Then there are two cases:
\\[1mm]
\begin{itemize}
\item[\textbf{Case 1:}]
$\zeta=(y,u(y))+\vartheta$ \;for a \,$y\in B_\varrho$ and a \,$\vartheta\in \omega(y)$ with
$|\vartheta|=\varepsilon$. \\ \\
We distinguish two subcases 1.\textit{i} and 1.\textit{ii:} \\[-1.3mm]
\begin{itemize}
\item[\textbf{1.\textit{i:}}]
The case $|x-y|\leq \frac{\cos \gamma}{L(1+\lambda+\cos \gamma)}$. \\ \\
As $u$ is $\lambda$-Lipschitz, we have \beqnn
|(x,u(x))-(y,u(y))|&\leq& (1+\lambda)|x-y| \\
&\leq& \frac{(1+\lambda)\cos \gamma}{L(1+\lambda+\cos \gamma)}. \eeqnn
Then \beqnn
|(x,u(x))-\zeta|&\geq& |\zeta-(y,u(y))|-|(x,u(x))-(y,u(y))| \\
&\geq& \frac{1}{L}\cos \gamma
-\frac{(1+\lambda)\cos \gamma}{L(1+\lambda+\cos \gamma)} \\
&=& \frac{\cos^2 \gamma}{L(1+\lambda+\cos \gamma)}. \eeqnn
\item[\textbf{1.\textit{ii:}}] The case
$|x-y|> \frac{\cos \gamma}{L(1+\lambda+\cos \gamma)}$. \\ \\
Again let $\pi^\bot((x,u(x)))$ be the orthogonal projection of $(x,u(x))$ onto
$(y,u(y))+\omega(y)$.
With (\ref{appendix1}) we estimate \beqn
|(x,u(x))-\zeta|&\geq& |(x,u(x))-\pi^\bot((x,u(x)))|\nonumber\\[2mm]
&\geq&|x-y|\cos \gamma \label{appendix3} \\
&>& \frac{\cos^2 \gamma}{L(1+\lambda+\cos \gamma)}.\nonumber \eeqn
\end{itemize} \vspace{3mm}
Both in Case 1.\textit{i} and Case 1.\textit{ii}  we have \beqn \label{case1angle}
|(x,u(x))-\zeta|\geq\frac{\cos^2 \gamma}{2L(1+\lambda)}. \eeqn \\
\item[\textbf{Case 2:}] $\zeta=(z,u(z))+\upsilon$
for a $z\in \partial B_\varrho$ and $\upsilon\in \omega(z)$
with $|\upsilon|\leq \varepsilon$.
\\ \\
As $x \in \overline{B}_{\frac{\varrho}{2}}$ we have $|x-z|\geq \frac{\varrho}{2}$.
Considering the orthogonal projection onto $(z,u(z))+\omega(z)$,
we estimate as in (\ref{appendix3}) \beqn \label{case2angle}
|(x,u(x))-\zeta| \geq \frac{\varrho}{2}\cos \gamma. \eeqn
\end{itemize} \vspace{5mm}
With (\ref{case1angle}) and (\ref{case2angle}) we have in any case \beqnn
|(x,u(x))-\zeta|\;\geq\;
\min\!\left\{ \frac{\varrho}{2}\cos \gamma, \,\frac{\cos^2 \gamma}{2L(1+\lambda)}\!\right\}.
\eeqnn
This proves part b) of Lemma \ref{tubularsigma}.
\end{itemize} \vspace{-7mm}
\hfill $\square$
\vspace{0.6cm}
\begin{lemma} \label{intersectionmostonepoint}
Let $f:M^m\rightarrow \R^{m+1}$ be an $(r,\lambda)$-immersion, $q\in M$ and $0<\varrho\leq r$.
Let $\omega\in G_{m+1,1}$ with $\R^{m+1}=\tau_f(p)\oplus \omega$ for all $p\in U_{\varrho,q}$.
Then for every $x\in \R^{m+1}$ the set $x+\omega$ intersects $f(U_{\varrho,q})$ in at most
one point. \end{lemma}
\textbf{Proof:} \\
After a rotation and a translation we may assume $f(U_{\varrho,q})=\{(y,u(y)):y\in B_\varrho\}$
with a $C^1$-function $u:B_\varrho\rightarrow \R$. Suppose the assertion of the lemma is false.
Then there is an $x\in \R^{m+1}$ such that $x+\omega$ intersects $f(U_{\varrho,q})$ in
$(y,u(y))$ and in $(z,u(z))$ with $y\neq z$. We may assume $y-z\in \R^1\times\{0\}\subset \R^m$.
With the same argument as in the paragraph after (\ref{anglelineest}) we conclude that there is a
$\xi\in \{(1-t)y+tz:\:t\in(0,1)\}\subset B_\varrho$ with
$\omega=\text{span}\{(e_1,\partial_1u(\xi))\}$. Moreover there is a unique $\zeta\in U_{\varrho,q}$
with $\tau_f(\zeta)=\text{span}\{(e_1,\partial_1u(\xi)),\ldots,(e_m,\partial_mu(\xi))\}$. Hence
$\omega\subset \tau_f(\zeta)$. But this contradicts
$\R^{m+1}=\tau_f(p)\oplus \omega$ for all $p\in U_{\varrho,q}$.
\hfill $\square$
\end{appendix}

\end{document}